\documentclass[a4paper, 11pt, twoside]{article}
\usepackage{latexsym,bm}
\usepackage[tbtags]{amsmath}
\usepackage{amssymb}

\usepackage{amsmath,amsfonts,amsthm,amssymb}

\usepackage{indentfirst}
\usepackage{paralist}
\usepackage[english]{babel}
\usepackage{amsmath,amssymb}
\usepackage{mathrsfs}

\usepackage{cite}
\usepackage{tabularx}
\usepackage{booktabs}
\usepackage{dcolumn}

\usepackage{multirow}

\usepackage{graphicx}

\usepackage{wrapfig}
\usepackage{multicol}
\usepackage{verbatim}

\newtheoremstyle{theorem}
  {10pt}
  {10pt}
  {\sl}
 {}
  {\bf}
  {. }
  { }
  {}
\theoremstyle{theorem}

\numberwithin{equation}{section}

\newtheoremstyle{defi}
  {10pt}
  {10pt}
  {\rm}
  {}
  {\bf}
  {. }
  { }
  {}
\theoremstyle{defi}

\allowdisplaybreaks \textheight 21.5 true cm
\usepackage{epsfig,graphicx,epsf,wrapfig}
\usepackage{dcolumn}
\usepackage{bm}
\begin{document}
\title{Stochastic 3D Navier-Stokes equations with nonlinear damping: martingale solution, strong solution and small time large deviation principles
\thanks
{The work is supported in part by a NSFC Grant No. 11531006,
PAPD of Jiangsu Higher Education Institutions and  National Basic Research Program of China (973 Program) No. 2013CB834100.
}}
\author{Hui Liu$^{1}$ ~~~~Hongjun Gao$^{1,2}$\thanks{Corresponding author, gaohj@njnu.edu.cn}}
\date{}
\maketitle

\noindent{\small{${^1}$Jiangsu Provincial Key Laboratory for NSLSCS, School of Mathematical Sciences, Nanjing Normal University, Nanjing 210023, PR China}}
\\
{\small{${^2}$ Institute of Mathematics, Jilin University, Changchun 130012, PR China}}
\\

\noindent \textbf{Abstract~}
In this paper, by using classical Faedo-Galerkin approximation and compactness method, the existence of martingale solutions for the stochastic 3D Navier-Stokes equations with nonlinear damping is obtained. The existence and uniqueness of strong solution are proved for $\beta > 3$ with any $\alpha>0$ and $\alpha \geq \frac12$ as $\beta = 3$. Meanwhile， a small time large deviation principle for the stochastic 3D Navier-Stokes equation with damping is proved for $\beta > 3$
with any $\alpha>0$ and $\alpha \geq \frac12$ as $\beta = 3$.
\\[2mm]
\textbf{Key words~}Stochastic Navier-Stokes equations, martingale solution, nonlinear damping, strong solution, large deviations.
\\[2mm]
\\
\textbf{2010 Mathematics Subject Classification~~~}34A37, 93B05, 93E03, 60H20, 34K50

\section{Introduction}
 We are concerned with very general classes of stochastic Navier-Stokes equation with damping. The damping is from the resistance to the motion of the flow. It describes various physical situations such as drag or friction effects, and some dissipative mechanisms\cite{bessaih1, hsiao}. Importance of such problems for climate modeling and physical fluid dynamics  is well known \cite{temam, griffies}.

 In this paper we consider the following stochastic three-dimensional Navier-Stokes equations with nonlinear damping
\begin{equation}\label{1.2}\begin{array}{l}
\left\{
\begin{array}{l}
 \frac{\partial u}{\partial t}-\mu\Delta u+(u\cdot\nabla)u+\nabla p+\alpha|u|^{\beta-1}u=f(t)+G(t,u)dW(t),\\
\nabla\cdot u=0,\\
u(t,x)|_{\partial D}=0,\\
u(x,0)=u_{0}(x).\\
\end{array}
\right.
\end{array}\end{equation}
Where $D\subset R^{3}$ be an open connected bounded subset with smooth boundary $\partial D$,  $u=(u_{1}, u_{2},u_{3})$ is the velocity,  $p$ is the pressure, $\beta>1$, $\alpha>0$ and $t\in[0,T]$, $u_{0}$ is the initial velocity, $f$ is the deterministic external force  and $G(t,u)dW(t)$, where $W$ is a cylindrical Wiener process, stands for the random forces.

The non-negative coefficients $\mu$ and $\alpha$ are called kinematic viscosity and sticky viscosity, respectively.
For a fixed $\alpha>0$, if $\mu>0$, these are called stochastic damped Navier-Stokes equations, whereas if $\mu=0$, they are the stochastic damped Euler equations.

The deterministic Navier-Stokes equation with nonlinear damping has been extensively investigated. For instance, Cai and Jiu have studied the the existence and regularity of solutions for three-dimensional Navier-Stokes equation with nonlinear damping \cite{cai}, they obtained the global weak solution for $\beta \ge 1$, the global strong solution for $\beta \ge \frac72$ and that the strong solution was unique for any $\frac72 \le \beta \le 5$,  based on it, Song and Hou considered the global attractor in \cite{song} and \cite{song1}. In \cite{jzd}, the authors considered the $L^2$ decay of weak solutions with $\beta \geq \frac{10}{3}$,  the optimal upper bounds of the higher-order derivative of the strong solution for $\frac72 \leq \beta < 5$ and  the asymptotic stability of the large solution to the system with $\beta \geq \frac72$. In \cite{wz}, the authors considered the regularity criterion of the 3D Navier-Stokes equations with nonlinear damping. In \cite{zhou}, for $\alpha = 1$, Zhou proved that the strong solution exists globally for $\beta \geq 3$ and strong-weak uniqueness for $\beta \geq 1$, and established two regularity criteria as $1\leq \beta \leq 3$. Oliveira has studied the existence of weak solutions for the generalized
Navier-Stokes equations with damping\cite{oliveira} for non-Newtonian fluids.

In a previous work, the existence of martingale solution of stochastic Euler equations was studied by many authors, for instance Capi$\acute{n}$ski and Cutland \cite{capinski} and  Bessaih \cite{bessaih2}, Brze$\acute{z}$niak and Peszat \cite{brz3}. Flandoli and  Gatarek \cite{flandoli1} obtained the martingale and stationary solutions for n-dimensional stochastic Navier-Stokes equations.  In the paper \cite{flandoli2, goldys}, the existence of martingale solutions and Markov selections of stochastic 3D Navier-Stokes equations were proved. An existence of martingale solution has been proved for the stochastic 3D Navier-Stokes equations with jump under proper assumptions \cite{dong}.

Large deviation have applications in many areas, such as in thermodynamics, stochastic mechanics, information theory and risk management, etc.,\cite{touchette}.
A small time large deviation principle for the two-dimensional stochastic Navier-Stokes equations driven by multiplicative noise was established in\cite{xu1}.  Using weak convergence method, a large deviation principle of Freidlin-Wentzell type for stochastic tamed 3D Navier-Stokes equations driven by multiplicative noise was proved in\cite{rockner2}. The small time large deviation principle for the stochastic 3D tamed Navier-Stokes equations was obtained in\cite{rockner3}.

Using Galerkin approximation, Prokhorov's theorem and Skorokhod's embedding theorem, the existence of martingale solutions for a stochastic 2D and 3D Navier-Stokes equations on unbounded domain is obtained\cite{brz2}. The presence of the nonlinear damping $|u|^{\beta-1}u$ in $(1.1)$ is beneficial to the regularity of the martingale solution. Using the nonlinear structure and delicated analysis we overcome some difficulties (such as $\int_{0}^{t}\int_{D}|u\cdot\nabla u|^{2}dxds$), we could prove the existence and uniqueness of strong solution for problem (1.1) for $\beta > 3$ for any $\alpha > 0$ and $\alpha \geq \frac12$ as $\beta = 3$. If we used the method in \cite{zhou} with further delicate analysis, we still need  the condition $\alpha \geq \frac12$ for $\beta \ge 3$, so our result improve the early results even for deterministic case.

As for the existence of local solution for stochastic 3D Navier-Stokes equations with damping,  we could use the monotonicity method as in \cite{liu-rockner}. Our purpose will prove the global existence of strong solution which depends on the various estimates of solution, the LDP and ergodicity etc also depend on such estimates of solution, so we use the Galerkin approximation, a priori estimate and weak convergence here.

The paper is organized as follows. In section 2, we recall basic definitions and introduce some auxiliary operators. Martingale solutions will be obtained in section 3. In section 4, we prove the existence and uniqueness of strong solution for $\beta > 3$ with any $\alpha>0$ and $\alpha \geq \frac12$ as $\beta = 3$. In section 5, we prove the small time large deviation principle for the stochastic 3D Navier-Stokes equation with damping for $\beta > 3$ with any $\alpha>0$ and $\alpha \geq \frac12$ as $\beta = 3$.

\section{Preliminaries}
In this section, we introduce some basic definitions, important lemmas which will be needed in this paper.

Let $D\subset R^{3}$ be an open subset with smooth boundary $\partial D$.
The norm in $L^{p}(D,R^{3})$ is defined by
\begin{equation*}\label{}\begin{array}{rl}
||u||_{L^{p}}=(\int_{D}|u(x)|^{p}dx)^{\frac{1}{p}},~~~~u\in L^{p}(D,R^{3}),~p>1.
\end{array}\end{equation*}
Let $H^{1}(D,R^{3})$ stands for the Sobolev space of all $u\in L^{2}(D,R^{3})$ for which there exists weak derivatives $\frac{\partial u}{\partial x_{i}}\in L^{2}(D,R^{3})$, $i=1,2,3$. It is a Hilbert space with the scalar product defined by
\begin{equation*}\label{}\begin{array}{rl}
(u,v)_{H^{1}}=( u, v)_{L^{2}}+((u,v)),~~u,~v\in H^{1}(D,R^{3}),
\end{array}\end{equation*}
where
\begin{equation*}\label{}\begin{array}{rl}
((u,v))=(\nabla u,\nabla v)_{L^{2}}=\sum\limits_{i=1}^{d}\int_{D}\frac{\partial u}{\partial x_{i}}\frac{\partial v}{\partial x_{i}}dx,~~u,~v\in H^{1}(D,R^{3}).
\end{array}\end{equation*}
Let $C^{\infty}_{c}(D,R^{3})$ be the space of all $R^{3}-$valued functions of class $\mathcal{C}^{\infty}$ with compact supports included in $D$ and let
\begin{equation*}\label{}\begin{array}{rl}
\mathcal{V}=\{u\in\mathcal{C}_{c}^{\infty}(D,R^{3}):divu=0\},\\
H=~the~closure~of~\mathcal{V}~in~L^{2}(D,R^{3}),\\
V=~the~ closure~of~ \mathcal{V} ~in ~H^{1}(D,R^{3}).
\end{array}\end{equation*}
For the equation $(1.1)$, we consider the abstract form as a stochastic evolution equation
\begin{equation}\label{}\begin{array}{rl}
du(t)+\mu \mathcal{A}u(t)dt+B(u(t))dt+g(u(t))dt=f(t)dt+G(u(t))dW(t)
\end{array}\end{equation}
with $u(0)=u_{0}$.\\
We define the bilinear operator $B(u,v)$, as $\langle B(u,v),z\rangle=b(u,v,z)=\int_{D}(u\cdot\nabla v)zdx$ for all $u,v,z\in V$, and also define the linear operator $\mathcal{A}$, as $\mathcal{A}u=-\Delta u$ and assume $\mu=1$ and $g(u)=\alpha|u|^{\beta-1}u$.\\
By the incompressibility condition, refer in Temam \cite{temam}
\begin{equation*}\label{}\begin{array}{rl}
b(u,w,v)=-b(u,v,w),~~~~u,v,w\in V.
\end{array}\end{equation*}
In particular,
\begin{equation*}\label{}\begin{array}{rl}
b(u,v,v)=0,~~u,v\in V.
\end{array}\end{equation*}
Let $(X,|\cdot|_{X})$, $(Y,|\cdot|_{Y})$ denotes two real normed spaces. The symbol $\mathcal{L}(X,Y)$ stands for the space of all bounded linear operators from $X$ to $Y$. The symbol $_{X'}\langle\cdot,\cdot\rangle_{X}$ denotes the standard duality pairing. If both spaces $X$ and $Y$ are separable Hilbert, $\mathcal{T}_{2}(Y,X)$ stands for the Hilbert space of all Hilbert-Schmidt operators from $Y$ to $X$ endowed with the standard norm.

For all $s>0$, the following standard scale of Hilbert spaces is defined by
\begin{eqnarray*}
V_{s}=the~closure~of~\mathcal{V}~in~H^{s}(D,R^{3}).
\end{eqnarray*}
Consider the natural embedding $j:V\hookrightarrow H$ and its adjoint $j^{*}:H\rightarrow V$. Suppose that $s>\frac{5}{2}$, we get that $V_{s}$ is dense in $V$ and the embedding $j_{s}:V_{s}\hookrightarrow V$ is continuous. By (Lemma 2.5\cite{holly}), then there exists a Hilbert space $U$ such that $U\subset V_{s}$, $U$ is dense in $V_{s}$ and the natural embedding $\iota_{s}:U\hookrightarrow V_{s}$ is compact.\\
We get
\begin{eqnarray*}
U\hookrightarrow V_{s}\hookrightarrow V\hookrightarrow H\cong H'\hookrightarrow V'\hookrightarrow V'_{s}\hookrightarrow U'.
\end{eqnarray*}
Fixing $n\in N$ and let $H_{n}=$ span$\{e_{1},e_{2},\cdots e_{n}\}$ denote the linear space spanned by the vectors $e_{1},e_{2},\cdots, e_{n}$. Let $P_{n}$ denote the operator from $U'$ to span$\{e_{1},e_{2},\cdots e_{n}\}$ given by
\begin{equation}\label{}\begin{array}{rl}
P_{n}u^{\ast}=\sum\limits^{n}_{i=1}{}_{U'}\langle u^{*},e_{i}\rangle_{U} e_{i},~~~u^{\ast}\in U'.
\end{array}\end{equation}

We consider the following functional spaces, similarly to those considered in \cite{mikulevicius}:

$\mathcal{C}([0,T];U'):=$ the space of continuous functions $u:[0,T]\rightarrow U'$ with the topology $\mathcal{T}_{1}$ induced by the norm $||u||_{\mathcal{C}([0,T],U')}=\sup_{t\in[0,T]}|u(t)|_{U'}$,

$L^{2}_{w}(0,T;V):=$ the space $L^{2}(0,T;V)$ with the weak topology $\mathcal{T}_{2}$,

$L^{2}(0,T;H):=$ the space of measurable functions $u:[0,T]\rightarrow H$ such that
\begin{eqnarray*}
q_{T}(u):=||u||_{L^{2}(0,T;H)}=(\int_{0}^{T}\int_{D}|u(t,x)|^{2}dxdt)^{\frac{1}{2}}<\infty,
\end{eqnarray*}
with the topology $\mathcal{T}_{3}$ generated by the seminorms $(q_{T})$.\\
Let $H_{w}$ be the Hilbert space $H$ endowed with the weak topology. Let

$\mathcal{C}([0,T];H_{w}):=$ the space of weakly continuous functions $u:[0,T]\rightarrow H$ endowed with the weakest topology $\mathcal{T}_{4}$ such that for any $h\in H$ the mappings $\mathcal{C}([0,T];H_{w})\ni u\rightarrow(u(\cdot),h)_{H}\in\mathcal{C}([0,T];R)$ are continuous.\\
We introduce the ball
\begin{align*}
\mathbb{B}=\{x\in H:|x|_{H}\leq r\}.
\end{align*}
Let $\mathbb{B}_{w}$ be the ball $\mathbb{B}$ endowed with the weak topology and $q$ be the metric compatible with the weak topology on $\mathbb{B}$. We will introduce the following subspace of the space $\mathcal{C}([0,T];H_{w})$
\begin{align*}
\mathcal{C}([0,T];\mathbb{B}_{w})=\{u\in\mathcal{C}([0,T];H_{w}):\sup\limits_{t\in[0,T]}|u(t)|_{H}\leq r\}.
\end{align*}
The space $\mathcal{C}([0,T];\mathbb{B}_{w})$ is metrizable with
\begin{align*}
\varrho(u,v)=\sup\limits_{t\in[0,T]}q(u(t),v(t)).
\end{align*}
{\bf Lemma 2.1}
(see~~Lemma~2.7~in~\cite{mikulevicius}) Let
$$\mathcal{Z}=\mathcal{C}([0,T];U')\cap L^{2}_{w}(0,T;V)\cap L^{2}(0,T;H)\cap \mathcal{C}([0,T];H_{w})$$
and let $\mathcal{T}$ denote the supremum of the corresponding topologies. Then a set $\mathcal{K}\subset \mathcal{Z}$ is $\mathcal{T}-$ relatively compact if the following three conditions hold

$(a)~~sup_{u\in \mathcal{K}}sup_{s\in[0,T]}|u(s)|_{H}<\infty,$

$(b)~~sup_{u\in \mathcal{K}}\int^{T}_{0}||u(s)||^{2}ds<\infty,$

$(c)~~lim_{\delta\rightarrow 0}sup_{u\in \mathcal{K}}sup_{{s,t\in[0,T]}}|u(t)-u(s)|_{U'}=0,~~~~\forall~|t-s|\leq\delta.$\\
Let $(\mathbb{S},\varrho)$ denotes a separable and complete metric space.\\
{\bf Definition 2.1}
Let $u\in\mathcal{C}([0,T],\mathbb{S})$. The modulus of continuity of $u$ on $[0,T]$ is given by
\begin{eqnarray*}
m(u,\delta)=\sup_{s,t\in[0,T],|t-s|\leq\delta}\varrho(u(t),u(s)),~~\delta>0.
\end{eqnarray*}
Let $(\Omega,\mathscr{F},\mathbb{P})$ denote a probability space with filtration $\mathbb{F}=(\mathscr{F})_{t\in[0,T]}$ satisfying the ordinary condition \cite{metivier1}, let $(X_{n})_{n\in N}$ denote a sequence of continous $\mathbb{F}-$adapted $\mathbb{S}-$valued processes.\\
{\bf Definition 2.2}
(\cite{metivier}) A sequence $(X_{n})_{n\in N}$ satisfies condition $[\mathbf{A}]$ iff $\forall~\varepsilon>0$, $\forall~\eta>0$,  $\exists~\delta>0$ such that for any sequence $(\tau_{n})_{n\in N}$ of $\mathbb{F}-$stopping times with $\tau_{n}\leq T$, we have
$$\sup\limits_{n\in \mathbb{N}} \sup\limits_{0\leq\theta\leq\delta}\mathbb{P}\{\varrho(X_{n}(\tau_{n}+\theta),X_{n}(\tau_{n}))\geq\eta\}\leq\varepsilon.$$
Applying the deterministic compactness criterion formulated in Lemma 2.1, we get the following Lemma 2.2.\\
{\bf Lemma 2.2}
(Tightness criterion~~\cite{metivier}) Assume that $(X_{n})_{n\in N}$ denote a sequence of continuous $\mathbb{F}-$adapted $U'-$valued processes such that

$(a)$ there exists a positive constant $C_{1}$ such that
$$\sup\limits_{n\in \mathbb{N}}\mathbb{E}[\sup\limits_{s\in[0,T]}|X_{n}(s)|^{2}_{H}]\leq C_{1},$$

$(b)$ there exists a positive constant $C_{2}$ such that
$$\sup\limits_{n\in \mathbb{N}}\mathbb{E}[\int^{T}_{0}||X_{n}(s)||^{2}ds]\leq C_{2},$$

$(c)$ $(X_{n})_{n\in N}$ satisfies the Aldous condition $[\mathbf{A}]$ in $U'$.\\
Let $\tilde{\mathbb{P}}_{n}$ be the law of $X_{n}$ on $\mathcal{Z}$. Then there exists a compact subset $K_{\varepsilon}$ of $\mathcal{Z}$ such that for any $\varepsilon>0$
$$\sup\limits_{n\in\mathbb{N}}\tilde{\mathbb{P}}_{n}(K_{\varepsilon})\geq1-\varepsilon.$$
{\bf Lemma 2.3}
(Lemma 2.4\cite{song}) Let $g(u)=\alpha|u|^{\beta-1}u$, then

(1) $g$ is continuously differential in $R^{3}$ and the Jacobian matrix is defined by, for any $u=(u_{1},u_{2},u_{3})$ in $R^{3}:$
\begin{equation*}\label{}\begin{array}{rl}
g'(u)=\alpha|u|^{\beta-3}\left(
                                                                    \begin{array}{ccc}
                                                                       (\beta-1)u_{1}^{2}+|u|^{2} & (\beta-1)u_{1}u_{2}  & (\beta-1)u_{1}u_{3} \\
                                                                       (\beta-1)u_{1}u_{2} & (\beta-1)u_{2}^{2}+|u|^{2} & (\beta-1)u_{2}u_{3} \\
                                                                       (\beta-1)u_{1}u_{3} & (\beta-1)u_{2}u_{3} & (\beta-1)u_{3}^{2}+|u|^{2} \\
                                                                     \end{array}
                                                                   \right
).
\end{array}\end{equation*}
Moreover, $g'(u)$ is positive definite and for any $u,v,w\in R^{3}:$
\begin{equation*}\label{}\begin{array}{rl}
|(g'(u)v)\cdot w|\leq c|u|^{\beta-1}|v||w|,
\end{array}\end{equation*}
where $c$ denote a positive constant depending on $\beta$ and $\alpha$.

(2) $g$ is monotonic in $R^{3}$, i.e., for any $u,v\in R^{3}:$
\begin{equation*}\label{}\begin{array}{rl}
(g(u)-g(v),u-v)\geq 0.
\end{array}\end{equation*}
{\bf Assumption A.} There exists positive constants $k_{1}$, $k_{2}$, $k_{3}$ such that
\begin{align*}
||\nabla G_{n}(u_{n}(r))P_{n}||^{2}_{\mathcal{T}_{2}(Y,H)}\leq k_{1}||u||^{2}+k_{2}|u|^{2}_{H}+k_{3}.
\end{align*}
\section{Existence of martingale solutions}
We shall consider the following stochastic evolution equation
\begin{equation}\label{1.2}\begin{array}{l}
\left\{
\begin{array}{l}
 du{(t)}+\mathcal{A}u(t)dt+B(u(t),u(t))dt+g(u)dt=f(t)dt+G(u(t))dW(t),\\
u(0)=u_{0}.\\
\end{array}
\right.
\end{array}\end{equation}
{\bf Definition 3.1}
We say that there exists a martingale solution of the equation $(3.1)$ iff there exists a stochastic basis $(\hat{\Omega},\hat{\mathscr{F}},\hat{\mathbb{F}},\hat{\mathbb{P}})$ with filtration $\hat{\mathbb{F}}=\{{\hat{\mathscr{F}}}_{t}\}_{t\geq0}$, a cylindrical Wiener process $\hat{W}$ on the space $Y$ and a progressively measurable process $u:[0,T]\times\hat{\Omega}\rightarrow H$, with $\hat{\mathbb{P}}-a.s.$ paths
\begin{equation*}\label{}\begin{array}{rl}
u(\cdot,\omega)\in C([0,T];H_{\omega})\cap L^{2}(0,T;V)\cap L^{\beta+1}(0,T; L^{\beta+1})
\end{array}\end{equation*}
such that $\hat{\mathbb{P}}-a.s.$, the identity
\begin{equation*}\label{}\begin{array}{rl}
(u(t),v)_{H}+\int^{t}_{0}\langle \mathcal{A}u(s),v\rangle ds+\int^{t}_{0}\langle B(u(s),u(s)),v\rangle ds+\int^{t}_{0}\langle g(u(s)),v\rangle ds
\end{array}\end{equation*}
\begin{equation}\label{}\begin{array}{rl}
=(u_{0},v)_{H}+\int^{t}_{0}\langle f(s),v\rangle ds+\langle\int^{t}_{0}G(u(s))d\hat{W}(s),v\rangle
\end{array}\end{equation}
holds true for any $t\in[0,T]$ and all $v\in\mathcal{V}$.\\
{\bf Theorem 3.1}
Let following assumptions $(A.1)-(A.3)$ be satisfied,

$(A.1)$ $W(t)$ be a cylindrical Wiener process in a separable Hilbert space $Y$ defined on the stochastic basis $(\Omega,\mathscr{F},\mathbb{F},\mathbb{P})$ with a filtration $\mathbb{F}=\{\mathscr{F}_{t}\}_{t\geq 0}$ and $H\subset Y$;

$(A.2)$ $u_{0}\in H,f\in L^{2}(0,T;V')$;

$(A.3)$ the mapping $G:V\rightarrow \mathcal{T}_{2}(Y,H)$ is Lipschitz continuous and
\begin{equation}\label{}\begin{array}{rl}
2\langle\mathcal{A}u,u\rangle-||G(u)||^{2}_{\mathcal{T}_{2}(Y,H)}\geq\eta||u||^{2}-\lambda_{0}|u|^{2}_{H}-\rho,
\end{array}\end{equation}
where, $u\in V$, for some constants $\lambda_{0}$, $\rho$ and $\eta\in (0,2]$. \\
Then there exists a martingale solution $(\hat{\Omega},\hat{\mathscr{F}},\hat{\mathbb{F}},\hat{\mathbb{P}},u)$ of problem $(3.1)$ such that
\begin{equation}\label{}\begin{array}{rl}
\hat{\mathbb{E}}[\sup\limits_{t\in[0,T]}|u(t)|^{2}_{H}+\int^{T}_{0}||u(t)||^{2}dt]< \infty.
\end{array}\end{equation}
However, $G$ extends to a mapping $G:H\rightarrow \mathcal{T}_{2}(Y,V')$ such that
\begin{equation}\label{}\begin{array}{rl}
||G(u)||^{2}_{\mathcal{T}_{2}(Y,V')}\leq C(1+|u|^{2}_{H}),~~u\in H.
\end{array}\end{equation}
Let $H_{n}=span{\{e_{i},1\leq i\leq n\}}$ be the orthogonal projection of $H$ on $H_{n}$, and let $P_{n}:U'\rightarrow H_{n}$ be defined by $(2.2)$.

Consider the following classical Faedo-Galerkin approximation in the space $H_{n}$ $:$
\begin{equation}\label{1.2}\begin{array}{l}
\left\{
\begin{array}{l}
 du_{n}{(t)}=-[P_{n}\mathcal{A}u_{n}(t)+B_{n}(u_{n}(t))+P_{n}g(u_{n})-P_{n}f(t)]dt+P_{n}G(u_{n}(t))dW(t),\\
u_{n}(0)=P_{n}u_{0}.\\
\end{array}
\right.
\end{array}\end{equation}
Using the It$\hat{o}$ formula and the Burkholder-Davis-Gundy inequality \cite{da1}, we will prove the following Lemma 3.1. Let $p$ satisfy the following condition
\begin{equation}\label{1.2}\begin{array}{l}
\left\{
\begin{array}{l}
 p\in[2,2+\frac{\eta}{2-\eta})~~if~\eta\in(0,2),\\
p\in[2,\infty)~~if~\eta=2.\\
\end{array}
\right.
\end{array}\end{equation}
Since $||u||:=||\nabla u||_{L^{2}}$ and $\langle\mathcal{A}u,u\rangle=((u,u))=(\nabla u,\nabla u)_{L^{2}}$,
\begin{equation*}\label{}\begin{array}{rl}
2\langle\mathcal{A}u,u\rangle-\eta||u||^{2}=(2-\eta)||u||^{2}
\end{array}\end{equation*}
for some constants $\eta\in(0,2]$, $\lambda_{0}\geq0$ and $\rho\in R$.

So inequality $(3.3)$ should be written equivalently in the following form
\begin{equation}\label{}\begin{array}{rl}
||G(u)||^{2}_{\mathcal{T}_{2}(Y,H)}\leq(2-\eta)||u||^{2}+\lambda_{0}|u|^{2}_{H}+\rho.
\end{array}\end{equation}
The following estimate and Lemma 3.2 are based on the similar method used by Brze$\acute{z}$niak and Motyl in \cite{brz2}.\\
{\bf Lemma 3.1}
There exist three positive constants $C_{1}(p)$, $C_{2}(p)$, $C_{3}$, for each $p\geq2$, such that
\begin{eqnarray}
\sup\limits_{n\geq1}\mathbb{E}(\sup\limits_{0\leq s\leq t}|u_{n}(s)|^{p}_{H})\leq C_{1}(p),
\end{eqnarray}
\begin{eqnarray}
\sup\limits_{n\geq1}\mathbb{E}\int^{t}_{0}|u_{n}(s)|^{p-2}_{H}||u_{n}(s)||^{2}ds\leq C_{2}(p),
\end{eqnarray}
\begin{eqnarray}
\sup\limits_{n\geq1}\mathbb{E}\int^{t}_{0}|u_{n}(s)|^{\beta+1}_{\beta+1}ds\leq C_{3}.
\end{eqnarray}

For any $n\in N$, the solution $u_{n}$ of the Galerkin equations defines a measure $\mathcal{L}(u_{n})$ on $(\mathcal{Z},\mathcal{T})$. Using Lemma 2.2 and (3.5), we will prove the tightness of this set of measures.\\
{\bf Lemma 3.2}
The family ${\mathcal{L}(u_{n})}$ is tight on $(\mathcal{Z},\mathcal{T})$.\\
Proof of Theorem 3.1.

By Lemma 3.2, we have the set of measures $\{\mathcal{L}(u_{n}),n\in N\}$ is tight on the space $(\mathcal{Z},\mathcal{T})$, and by [Corollary 3.12\cite{brz2}], then there exists a subsequence $(n_{k})_{k}$, a probability space $(\tilde{\Omega},\tilde{\mathscr{F}},\tilde{\mathbb{P}})$, and $\mathcal{Z}-$valued random variables $\tilde{u}$, $\tilde{u}_{n_{k}}$, $k\geq1$ such that

~~~~~~~~$\tilde{u}_{n_{k}}$ has the same law as $u_{n_{k}}$ on $\mathcal{Z}$ and $\tilde{u}_{n_{k}}\rightarrow u_{n_{k}}$ in $\mathcal{Z}$,  $\tilde{\mathbb{P}}-$a.s.\\
We denote the subsequence $(\tilde{u}_{n_{k}})_{k}$ by $(\tilde{u}_{n})_{n}$.

Let us define the process $\tilde{M}_{n}(t)$ with trajectories in $\mathcal{C}([0,T];H)$
\begin{equation*}\label{}\begin{array}{rl}
\tilde{M}_{n}(t)=\tilde{u}_{n}(t)-P_{n}\tilde{u}_{0}+\int^{t}_{0}P_{n}\mathcal{A}\tilde{u}_{n}(s)ds+\int^{t}_{0}B_{n}(\tilde{u}_{n}(s))ds
\end{array}\end{equation*}
\begin{equation}\label{}\begin{array}{rl}
+\int^{t}_{0}P_{n}g(\tilde{u}_{n}(s))ds-\int^{t}_{0}P_{n}f(s)ds,~~t\in[0,T],~~n\geq1.
\end{array}\end{equation}

The following proof is based on the similar method used by Brze$\acute{z}$niak and Motyl in \cite{brz2}. Since $\tilde{u}_{n}$ and $u_{n}$ have the same laws, for all $s,t\in[0,T]$, $s\leq t$ all functions $h$ bounded continuous on $\mathcal{C}([0,s];U')$, and all $\psi,\zeta\in U$.\\
{\bf Lemma 3.3}
(Lemma 5.5\cite{brz2}) For all $s,t\in[0,T]$ such that $s\leq t$ and all $\psi\in U$:

(a)$\lim_{n\rightarrow\infty}(\tilde{u}_{n},P_{n}\psi)_{H}=(\tilde{u},\psi)_{H},~\tilde{\mathbb{P}}-a.s.,$

(b)$\lim_{n\rightarrow\infty}\int_{s}^{t}\langle\mathcal{A}\tilde{u}_{n}(\sigma),P_{n}\psi\rangle d\sigma=\int_{s}^{t}\langle\mathcal{A}\tilde{u}(\sigma),\psi\rangle d\sigma,~\tilde{\mathbb{P}}-a.s.,$

(c)$\lim_{n\rightarrow\infty}\int_{s}^{t}\langle B(\tilde{u}_{n}(\sigma)),P_{n}\psi\rangle d\sigma=\int_{s}^{t}\langle B(\tilde{u}(\sigma)),\psi\rangle d\sigma,~\tilde{\mathbb{P}}-a.s..$\\
{\bf Lemma 3.4} For all $s,t\in[0,T]$ such that $s\leq t$ and all $\psi\in U$:
\begin{eqnarray*}
(a)~~lim_{n\rightarrow\infty}(g(\tilde{u}_{n}(t)),P_{n}\psi)_{H}=(g(\tilde{u}(t)),\psi)_{H},~~~~\tilde{\mathbb{P}}-a.s.
\end{eqnarray*}
Proof. Let us fix $s,t\in[0,T]$, $s\leq t$ and $\psi\in U$. By Lemma 3.3(a), we get that $\tilde{u}_{n}\rightarrow \tilde{u}$ in $C([0,T],H_{w})$, $\tilde{\mathbb{P}}-$a.s. and by (2.2), we get $P_{n}\psi\rightarrow\psi$ in $H$.\\
And because
\begin{align*}
|g(\tilde{u}_{n}(t))-g(\tilde{u}(t))|&=|\alpha|\tilde{u}_{n}|^{\beta-1}\tilde{u}_{n}-\alpha|\tilde{u}|^{\beta-1}\tilde{u}|\\
&\leq\alpha|\tilde{u}_{n}|^{\beta-1}|\tilde{u}_{n}-\tilde{u}|+\alpha||\tilde{u}_{n}|^{\beta-1}-|\tilde{u}|^{\beta-1}||\tilde{u}|\\
&\leq\alpha|\tilde{u}_{n}|^{\beta-1}|\tilde{u}_{n}-\tilde{u}|+C\alpha|\tilde{u}|(|\tilde{u}_{n}|^{\beta-2}+|\tilde{u}|^{\beta-2})|\tilde{u}_{n}-\tilde{u}|.
\end{align*}
We get that $g(\tilde{u}_{n}(t))\rightarrow g(\tilde{u}(t))$. We infer that assertion (a) holds.\\
{\bf Lemma 3.5}
(Lemma 5.6\cite{brz2}) For all $s,t\in[0,T]$ such that $s\leq t$ and all $\psi\in U$:
$$\lim_{n\rightarrow\infty}\mathbb{E}[\langle \tilde{M}_{n}(t)-\tilde{M}_{n}(s),\psi\rangle h(\tilde{u}_{n|[0,s]})]=\mathbb{E}[\langle \tilde{M}(t)-\tilde{M}(s),\psi\rangle h(\tilde{u}_{|[0,s]})].$$
{\bf Lemma 3.6}
(Lemma 5.7\cite{brz2}) For all $s,t\in[0,T]$ such that $s\leq t$ and all $\psi,\zeta\in U$:
\begin{equation*}\label{}\begin{array}{rl}
\lim_{n\rightarrow\infty}\mathbb{E}[(\langle \tilde{M}_{n}(t),\psi\rangle\langle \tilde{M}_{n}(t),\zeta\rangle-\langle \tilde{M}_{n}(s),\psi\rangle\langle \tilde{M}_{n}(s),\zeta\rangle)h(\tilde{u}_{n|[0,s]})]\\
=\mathbb{E}[(\langle \tilde{M}(t),\psi\rangle\langle \tilde{M}(t),\zeta\rangle-\langle \tilde{M}(s),\psi\rangle\langle \tilde{M}(s),\zeta\rangle)h(\tilde{u}_{|[0,s]})].
\end{array}\end{equation*}
{\bf Lemma 3.7}
(Lemma 5.8\cite{brz2}) For all $s,t\in[0,T]$ such that $s\leq t$ and all $\psi,\zeta\in U$, we have
\begin{equation*}\label{}\begin{array}{rl}
\lim_{n\rightarrow\infty}\mathbb{E}[(\int^{t}_{s}(G(\tilde{u}_{n}(\sigma))^{\ast}P_{n}\psi,G(\tilde{u}_{n}(\sigma))^{\ast}P_{n}\zeta)_{Y}d\sigma)\cdot h(\tilde{u}_{n|[0,s]})]\\
=\mathbb{E}[(\int^{t}_{s}(G(\tilde{u}(\sigma))^{\ast}\psi,G(\tilde{u}(\sigma))^{\ast}\zeta)_{Y}d\sigma)\cdot h(\tilde{u}_{|[0,s]})].
\end{array}\end{equation*}
Now, we apply the idea analogous to which used by Da Prato and Zabczyk (\cite{da2}, Section 8.3),  \cite{holly} and \cite{brz2}.
By the martingale representation theorem\cite{da2}, we will get conclude. This completes the proof of Theorem 3.1.
\section{Existence and uniqueness of strong solution}
In order to get the strong solution for the stochastic Navier-Stokes equation with damping, we have the following crucial lemma.

We first introduce the Galerkin systems associated to the original equation and establish some uniform a priori estimates. Recall that
\begin{align*}
\{e_{1},e_{2},\cdots\}\subset V
\end{align*}
is an orthonormal basis of $H$.
For any $n\geq1$, let $H_{n}=span\{e_{1},e_{2},\cdots,e_{n}\}$ and let $P_{n}:V'\rightarrow H_{n}$ be defined by
\begin{align}
P_{n}u=\sum\limits_{i=1}^{n} {}_{V'}\langle u,e_{i}\rangle_{V}e_{i},~~u\in V'.
\end{align}
Hence, $P_{n}|_{H}$ denote the orthogonal projection onto $H_{n}$ in $H$.

We consider the classical Faedo-Galerkin approximation in the space $H_{n}$ defined by
\begin{align}
 du_{n}{(t)}=-[P_{n}\mathcal{A}u_{n}(t)+B_{n}(u_{n}(t))+P_{n}g(u_{n})-P_{n}f(t)]dt+P_{n}G(u_{n}(t))dW(t)
\end{align}
for $u_{n}(0)=P_{n}u(0)$.\\
The following lemma is the main preliminary step in the proof of Theorem 4.2.\\
{\bf Lemma 4.1} Suppose that $\beta > 3$ with any $\alpha>0$ and $\alpha \geq \frac12$ as $\beta = 3$, $u_{n}(0)\in V\cap L^{\beta+1}(D)$, and $\nabla f\in L^{2}(\Omega;L^{2}(0,T;H))$, then there exists  a positive constant $C$, such that
 \begin{align*}
\sup\limits_{n\geq1}\mathbb{E}(\sup_{0\leq t\leq T}|\nabla u_{n}(t)|^{2}_{2}&+\int_{0}^{T}||\nabla u_{n}(r)||^{2}dr+ \int_{0}^{T}\int_{D}|u_{n}|^{\beta-1}|\nabla u_{n}|^{2}dsdr\notag\\
&+\int_{0}^{T}
|\nabla |u_{n}|^{\frac{\beta+1}{2}}|^{2}_{2}dr)\leq C(\mathbb{E}|\nabla u_{n}(0)|^{2}_{2}+1).
\end{align*}
Proof. We apply It$\hat{o}$ formula to $|\nabla u_{n}|^{2}_{2}$ for $t\in[0,T]$,
\begin{align}
&|\nabla u_{n}(t)|^{2}_{2}+2\int_{0}^{t}||\nabla u_{n}(r)||^{2}dr\notag\\
&+2\alpha\int_{0}^{t}\int_{D}|u_{n}|^{\beta-1}|\nabla u_{n}|^{2}dsdr+\frac{8\alpha(\beta-1)}{(\beta+1)^{2}}\int_{0}^{t}
|\nabla |u_{n}|^{\frac{\beta+1}{2}}|^{2}_{2}dr\notag\\
&\leq |P_{n}\nabla u_{n}(0)|^{2}_{2}+2\int_{0}^{t}|\langle\Delta u_{n}, B(u_{n})\rangle|dr+2\int_{0}^{t}|\langle\Delta u_{n},f\rangle|dr\notag\\
&+2|\int_{0}^{t}\langle G_{n}(u_{n}(r))dW(r),\Delta u_{n}(r)\rangle|+\int_{0}^{t}||\nabla G_{n}(u_{n}(r))P_{n}||^{2}_{\mathcal{T}_{2}(Y,H)}dr\notag\\
&=|P_{n}\nabla u_{n}(0)|^{2}_{2}+\sum\limits_{i=1}^{4}J_{i}(t).
\end{align}
Since
\begin{align}
0<\frac{2}{\beta-1}<1,~~~~~for~~\beta>3,
\end{align}
using Young's inequality and (4.4) to estimate $J_{1}(t)$, we deduce
\begin{align}
J_{1}(t)&\leq2\int_{0}^{t}(\frac{1}{2}|\Delta u_{n}|^{2}_{2}+\frac{1}{2}|u_{n}\cdot\nabla u_{n}|^{2}_{2})dr\notag\\
&\leq \int_{0}^{t}|\Delta u_{n}|^{2}_{2}dr+\int_{0}^{t}\int_{D}|\nabla u_{n}|^{2}|u_{n}|^{2}dsdr\notag\\
&=\int_{0}^{t}||\nabla u_{n}||^{2}dr+\int_{0}^{t}\int_{D}|u_{n}|^{2}|\nabla u_{n}|^{\frac{4}{\beta-1}}|\nabla u_{n}|^{2-\frac{4}{\beta-1}}dsdr\notag\\
&\leq \int_{0}^{t}||\nabla u_{n}||^{2}dr+\int_{0}^{t}[\int_{D}(|u_{n}|^{2}|\nabla u_{n}|^{\frac{4}{\beta-1}})^{\frac{\beta-1}{2}}ds]^{\frac{2}{\beta-1}}[\int_{D}(|\nabla u_{n}|^{2-\frac{4}{\beta-1}})^{\frac{\beta-1}{\beta-3}}ds]^{\frac{\beta-3}{\beta-1}}dr\notag\\
&\leq \int_{0}^{t}||\nabla u_{n}||^{2}dr+\varepsilon\int_{0}^{t}\int_{D}|\nabla u_{n}|^{2}|u_{n}|^{\beta-1}dsdr+C(\varepsilon)\int_{0}^{t}|\nabla u_{n}(r)|^{2}_{2}dr.
\end{align}
The H$\ddot{o}$lder's inequality and Young's inequality imply that
\begin{align}
J_{2}(t)\leq\frac{1}{4}\int_{0}^{t}||\nabla u_{n}(r)||^{2}dr+4(|\nabla f|^{2}_{L^{2}(0,t;H)}).
\end{align}
Taking the supremum and expectation over the interval $[0,t]$ on the equality (4.3), we estimate the last two items.
By the Assumption A, we have
\begin{align*}
||\nabla G_{n}(u_{n})||^{2}_{\mathcal{T}_{2}(Y,H)}\leq k_{1}||u_{n}||^{2}+k_{2}|u_{n}|^{2}_{H}+k_{3}.
\end{align*}
Now applying Burkholder-Davies-Gundy inequality and Young's inequality, we have
\begin{align}
&\mathbb{E}(\sup_{0\leq s\leq t}|J_{3}(s)|)\leq2\sqrt{2}\mathbb{E}\{\int_{0}^{t}|\nabla u_{n}(r)|^{2}_{2}|\nabla G_{n}(u_{n}(r))P_{n}|^{2}_{\mathcal{T}_{2}(Y,H)}dr\}^{\frac{1}{2}}\notag\\
&\leq\frac{1}{2}\mathbb{E}\sup_{0\leq s\leq t}|\nabla u_{n}|^{2}_{2}+4\mathbb{E}\int_{0}^{t}||\nabla G_{n}(u_{n}(r))||^{2}_{\mathcal{T}_{2}(Y,H)}dr\notag\\
&\leq\frac{1}{2}\mathbb{E}\sup_{0\leq s\leq t}|\nabla u_{n}|^{2}_{2}+4\mathbb{E}\int_{0}^{t}(k_{1}||u_{n}||^{2}+k_{2}|u_{n}|^{2}_{H}+k_{3} )dr.
\end{align}
Similarly, we have
\begin{align}
\mathbb{E}(\sup_{0\leq s\leq t}|J_{4}(s)|)\leq \mathbb{E}\int_{0}^{t}(k_{1}||u_{n}||^{2}+k_{2}|u_{n}|^{2}_{H}+k_{3} )dr.
\end{align}
Substituting (4.5)-(4.8) into (4.3), we have
\begin{align}
&\mathbb{E}\sup_{0\leq s\leq t}|\nabla u_{n}(s)|^{2}_{2}+\frac{3}{4}\mathbb{E}\int_{0}^{t}||\nabla u_{n}(r)||^{2}dr\notag\\
&+2\alpha \mathbb{E}\int_{0}^{t}\int_{D}|u_{n}|^{\beta-1}|\nabla u_{n}|^{2}dsdr+\frac{8\alpha(\beta-1)}{(\beta+1)^{2}}\mathbb{E}\int_{0}^{t}
|\nabla |u_{n}|^{\frac{\beta+1}{2}}|^{2}_{2}dr\notag\\
&\leq \mathbb{E}|P_{n}\nabla u_{n}(0)|^{2}_{2}+\varepsilon\mathbb{E}\int_{0}^{t}\int_{D}|\nabla u_{n}|^{2}|u_{n}|^{\beta-1}dsdr+C(\varepsilon)\mathbb{E}\int_{0}^{t}|\nabla u_{n}(r)|^{2}_{2}dr\notag\\
&+\frac{1}{2}\mathbb{E}\sup_{0\leq s\leq t}|\nabla u_{n}|^{2}_{2}+5\mathbb{E}\int_{0}^{t}(k_{1}||u_{n}||^{2}+k_{2}|u_{n}|^{2}_{H}+k_{3})dr+C|\nabla f|^{2}_{L^{2}(0,t;H)}.
\end{align}
We have
\begin{align}
&\mathbb{E}\sup_{0\leq s\leq t}|\nabla u_{n}(s)|^{2}_{2}+\frac{3}{2}\mathbb{E}\int_{0}^{t}||\nabla u_{n}(r)||^{2}dr\notag\\
&+2(2\alpha-\varepsilon) \mathbb{E}\int_{0}^{t}\int_{D}|u_{n}|^{\beta-1}|\nabla u_{n}|^{2}dsdr+\frac{16\alpha(\beta-1)}{(\beta+1)^{2}}\mathbb{E}\int_{0}^{t}
|\nabla |u_{n}|^{\frac{\beta+1}{2}}|^{2}_{2}dr\notag\\
&\leq 2\mathbb{E}|P_{n}\nabla u_{n}(0)|^{2}_{2}+C\mathbb{E}\int_{0}^{t}|\nabla u_{n}(r)|^{2}_{2}dr\notag\\
&+10\mathbb{E}\int_{0}^{t}(k_{1}||u_{n}||^{2}+k_{2}|u_{n}|^{2}_{H}+k_{3})dr+C|\nabla f|^{2}_{L^{2}(0,t;H)}\notag\\
&\leq 2\mathbb{E}|P_{n}\nabla u_{n}(0)|^{2}_{2}+C\int_{0}^{t}\mathbb{E}\sup\limits_{0\leq\tau\leq r}|\nabla u_{n}(\tau)|^{2}_{2}dr+C(t)(1+|\nabla f|^{2}_{L^{2}(0,t;H)}).
\end{align}
Choosing sufficiently small $\varepsilon$ and Gronwall's inequality, we deduce
\begin{align}
\sup\limits_{n\geq1}\mathbb{E}(\sup_{0\leq t\leq T}|\nabla u_{n}(t)|^{2}_{2}&+\int_{0}^{T}||\nabla u_{n}(r)||^{2}dr+ \int_{0}^{T}\int_{D}|u_{n}|^{\beta-1}|\nabla u_{n}|^{2}dsdr\notag\\
&+\int_{0}^{T}
|\nabla |u_{n}|^{\frac{\beta+1}{2}}|^{2}_{2}dr)\leq C(\mathbb{E}|\nabla u_{n}(0)|^{2}_{2}+1).
\end{align}
For $\alpha \geq \frac12$ as $\beta = 3$, we can get the above estimate easily. This completes the proof of Lemma 4.1.

Now we will study the weak convergence of approximating sequences. Using the priori estimate in the Lemma 3.1 and Lemma 4.1, we deduce the
existence of a subsequence of Galerkin elements $u_{n}$ and process $u\in L^{2}([0,T]\times\Omega,V)\cap L^{\beta+1}([0,T]\times\Omega,L^{\beta+1})\cap
L^{2}(\Omega,L^{\infty}([0,T],H))$, $\nabla u\in L^{2}([0,T]\times\Omega,V)\cap L^{2}(\Omega,L^{\infty}([0,T],H))$, $-\mathcal{A}u-B(u)+f\in L^{2}([0,T]\times\Omega,V')$, $g\in L^{s}([0,T]\times\Omega,L^{s}(D)),~~s=(\beta+1)^{*}=\frac{\beta+1}{\beta}$ and $G\in L^{2}([0,T]\times\Omega,\mathcal{T}_{2})$ for which the following limits hold:

(i) $u_{n}\rightarrow u$ weakly in $L^{2}([0,T]\times\Omega,V)$, $\nabla u_{n}\rightarrow\nabla u$ weakly in $L^{2}([0,T]\times\Omega,V)$,

(ii) $u_{n}\rightarrow u$ weakly in $L^{\beta+1}([0,T]\times\Omega,L^{\beta+1})$,

(iii) $u_{n}$ is weak star converging to $u$ in $L^{2}(\Omega,L^{\infty}([0,T],H))$, $\nabla u_{n}$ is weak star converging to $\nabla u$ in $L^{2}(\Omega,L^{\infty}([0,T],H))$,

(iv) $-\mathcal{A}u_{n}-B(u_{n})+f\rightarrow-\mathcal{A}u-B(u)+f$ weakly in $L^{2}([0,T]\times\Omega,V')$,

(v) $G_{n}(u_{n})P_{n}\rightarrow G$ weakly in $L^{2}([0,T]\times\Omega,\mathcal{T}_{2})$,

(vi) $P_{n}g(u_{n})\rightarrow g$ weakly in $g\in L^{s}([0,T]\times\Omega,L^{s}(D))$.
\subsection{Uniqueness of strong solution}
Let us recall that by assumption (A.3) the mapping $G:V\rightarrow \mathcal{T}_{2}(Y,H)$ is Lipschitz continuous, i.e. for some $L>0$ the following inequality holds
\begin{equation}\label{}\begin{array}{rl}
||G(u_{1}(s))-G(u_{2}(s))||^{2}_{\mathcal{T}_{2}(Y,H)}\leq L||u_{1}(s)-u_{2}(s)||,~~s\in[0,T].
\end{array}\end{equation}
{\bf Definition 4.1 (Pathwise uniqueness)} We say that the pathwise uniqueness holds for Eq. (3.1) if whenever we are given two weak solutions of Eq. (3.1) defined on the same probability space together with the same Brownian motion
\begin{equation*}\label{}\begin{array}{rl}
(\Omega,\mathscr{F},P; (\mathscr{F})_{t\geq0};W;u_{1}),
\end{array}\end{equation*}
\begin{equation*}\label{}\begin{array}{rl}
(\Omega,\mathscr{F},P; (\mathscr{F})_{t\geq0};W;u_{2}),
\end{array}\end{equation*}
the condition $P\{u_{1}(0)=u_{2}(0)\}=1$ implies $P\{\omega:u_{1}(t,\omega)=u_{2}(t,\omega),\forall t\geq0\}=1$.\\
The following Yamada-Watanabe theorem holds in this case(\cite{rockner, rockner1}).\\
{\bf Theorem 4.1} Under the conditions of Lemma 4.1,  the existence of martingale solutions plus pathwise uniqueness implies the existence of a unique strong solution.\\
{\bf Theorem 4.2} Under the conditions of Lemma 4.1, and suppose that assumption (A.1)-(A.3) are satisfied with  $\beta > 3$ with any $\alpha>0$ and $\alpha \geq \frac12$ as $\beta = 3$ and  $L<2$, then pathwise uniqueness holds for (3.1).\\
Proof. Let $u_{1}$ and $u_{2}$ be two solutions of Eq. (3.1) defined on the same probability space together with the same Brownian motion and starting from the same initial value. \\
Set
\begin{equation*}\label{}\begin{array}{rl}
U=u_{1}-u_{2}.
\end{array}\end{equation*}
Then $U$ satisfies the following equation
\begin{equation*}\label{}\begin{array}{rl}
dU(t)+[\mathcal{A}U(t)+(B(u_{1}(t))-B(u_{2}(t)))+(g(u_{1}(t))-g(u_{2}(t)))]dt
\end{array}\end{equation*}
\begin{equation*}\label{}\begin{array}{rl}
=[G(u_{1}(t))-G(u_{2}(t))]dW(t).
\end{array}\end{equation*}
Let $r(t)=a\int_{0}^{t}||u_{2}||^{4}ds$, $t\in[0,T]$, where $a$ is a positive constant. Let
$$F(t,\chi)=e^{-r(t)}|\chi|_{H}^{2},~~~~(t,\chi)\in[0,T]\times H.$$
Then by It$\hat{o}$'s formula, we have
\begin{align}
e^{-r(t)}|U(t)|_{H}^{2}&=\int_{0}^{t}e^{-r(s)}\{-r'(s)|U(s)|_{H}^{2}\notag\\
&-2\langle\mathcal{A}U(s)+(B(u_{1}(s))-B(u_{2}(s)))+(g(u_{1}(s))-g(u_{2}(s))),U(s)\rangle \}ds\notag\\
&+\frac{1}{2}\int_{0}^{t}Tr[(G(u_{1}(s))-G(u_{2}(s)))\frac{\partial^{2}F}{\partial \chi^{2}}(G(u_{1}(s))-G(u_{2}(s)))^{*}]ds\notag\\
&+2\int_{0}^{t}e^{-r(s)}\langle G(u_{1}(s))-G(u_{2}(s)), U(s)dW(s)\rangle\notag\\
&\leq\int_{0}^{t}e^{-r(s)}[-r'(s)|U(s)|_{H}^{2}-2||U(s)||^{2}\notag\\
&-2\langle B(u_{1}(s))-B(u_{2}(s)),U(s)\rangle-2\langle g(u_{1}(s))-g(u_{2}(s)),U(s)\rangle]ds\notag\\
&+\int_{0}^{t}e^{-r(s)}|| G(u_{1}(s))-G(u_{2}(s))||^{2}_{\mathcal{T}_{2}(Y,H)}ds\notag\\
&+2\int_{0}^{t}e^{-r(s)}\langle G(u_{1}(s))-G(u_{2}(s)), U(s)dW(s)\rangle.
\end{align}
We have
\begin{eqnarray*}
B(u_{1}(s))-B(u_{2}(s))=B(u_{1}(s),U(s))+B(U(s),u_{2}(s)),~~s\in[0,T].
\end{eqnarray*}
So
\begin{eqnarray*}
\langle B(u_{1}(s))-B(u_{2}(s)),U(s)\rangle=\langle B(U(s),u_{2}(s)),U(s)\rangle,
\end{eqnarray*}
and hence
\begin{eqnarray*}
|2\langle B(u_{1}(s))-B(u_{2}(s)),U(s)\rangle|\leq C |U|^{\frac{1}{2}}_{H}||U||^{\frac{3}{2}}||u_{2}||,~~s\in[0,T].
\end{eqnarray*}
Therefore for every $\varepsilon>0$ there exists $C_{\varepsilon}>0$ such that
\begin{eqnarray}
|2\langle B(u_{1}(s))-B(u_{2}(s)),U(s)\rangle|\leq \varepsilon||U||^{2}+C_{\varepsilon}|U|^{2}_{H}||u_{2}||^{4},~~s\in[0,T].
\end{eqnarray}
By Lemma (2.3), we have
\begin{align}
2\langle g(u_{1}(s))-g(u_{2}(s)),U(s)\rangle\geq0.
\end{align}
Let $a=C_{\varepsilon}$, we get
\begin{align*}
-r'(s)|U(s)|_{H}^{2}+C_{\varepsilon}|U|^{2}_{H}||u_{2}||^{4}=-a|U|^{2}_{H}||u_{2}||^{4}+C_{\varepsilon}|U|^{2}_{H}||u_{2}||^{4}=0.
\end{align*}
By (4.12), (4.13), (4.14) and (4.15), we have
\begin{align}
&e^{-r(t)}|U(t)|_{H}^{2}\notag\\
&+\int_{0}^{t}e^{-r(s)}(2-\varepsilon-L)||U(s)||^{2}ds\notag\\
&\leq 2\int_{0}^{t}e^{-r(s)}\langle G(u_{1}(s))-G(u_{2}(s)), U(s)dW(s)\rangle.
\end{align}
Let us choose $\varepsilon>0$ such that $2-\varepsilon-L>0$, we have
\begin{align}
e^{-r(t)}|U(t)|_{H}^{2}\leq2\int_{0}^{t}e^{-r(s)}\langle G(u_{1}(s))-G(u_{2}(s)), U(s)dW(s)\rangle.
\end{align}
Taking expectations for (4.17), we get
\begin{align}
\mathbb{E}(e^{-r(t)}|U(t)|_{H}^{2})\leq0.
\end{align}
By Lemma 4.1, we have $\mathbb{E}[\int_{0}^{T}||u_{2}||^{4}dt]<\infty$ for $\beta > 3$ with any $\alpha>0$ and $\alpha \geq \frac12$ as $\beta = 3$, and (4.18), we get
\begin{align}
\mathbb{E}(|U(t)|_{H}^{2})\leq0.
\end{align}
This completes the proof of Theorem 4.2.\\
{\bf Remark} If $W(\cdot)$ be a $H-$valued Brownian motion on $(\Omega,\mathscr{F},P)$ with the covariance operator $Q$, which is a positive, symmetric, trace class operator on $H$ and $G(t,u)$ satisfy the following hypotheses (A.1)-(A.4), using similar method, we get existence and uniqueness of strong solution on the space
$u\in L^{2}(\Omega;C([0,T];H))\cap L^{2}(\Omega\times[0,T];V)\cap L^{\beta+1}(\Omega\times[0,T];L^{\beta+1})$ and $\nabla u\in L^{2}([0,T]\times\Omega;V)\cap L^{2}(\Omega,C([0,T];H))$ for $\beta > 3$ with any $\alpha>0$ and $\alpha \geq \frac12$ as $\beta = 3$.
\section{Large deviation principle}\label{s1}
Let $W(\cdot)$ be a $H-$valued Brownian motion on $(\Omega,\mathscr{F},P)$ with the covariance operator $Q$, which is a positive, symmetric, trace class operator on $H$. Assume that $L_{Q}(H;H)$ denote the space of linear operator $S$ such that $SQ^{\frac{1}{2}}$ is a Hilbert-Schmidt operator from $H$ to $H$. Define the norm on the space $L_{Q}(H;H)$ by $|S|_{L_{Q}}=\sqrt{Tr(SQS^{*})}$. Assume that $L_{Q}^{V}$ denote the class of linear operators $\tilde{S}$ such that $\tilde{S}Q^{\frac{1}{2}}$ is a Hilbert-Schmidt operator from $H$ to $V$, endowed with the norm $|\tilde{S}|^{2}_{L_{Q}^{V}}=Tr(\tilde{S}Q\tilde{S}^{*})$.\\
In this section, we denote the following inner products
$$(u,v)=\int_{D}u\cdot vdx,~~\forall u,v\in H,$$
$$((u,v))=\sum\limits_{i=1}^{3}\int_{D}\nabla u_{i}\cdot\nabla v_{i}dx,~~\forall u,v\in V,$$
and norms $|\cdot|_{2}=(\cdot,\cdot)^{\frac{1}{2}}$, $||\cdot||=((\cdot,\cdot))$ and let $f(t)=0$.\\
Consider the following hypotheses\\
(A.1) $\mathbb{E}(|\xi|^{4})<\infty$.\\
(A.2) There exists a constant $L$ such that $|G(t,u)|^{2}_{L_{Q}}\leq L(1+|u|^{2})$, for all $t\in(0,T)$, and all $u\in H$.\\
(A.3) There exists a constant $\hat{L}$ such that $|G(t,u)|^{2}_{L_{Q}^{V}}\leq \hat{L}(1+||u||^{2})$, for all $t\in(0,T)$, and all $u\in V$.\\
(A.4) There exists a constant $K$ such that $|G(t,u)-G(t,v)|^{2}_{L_{Q}}\leq K|u-v|^{2}$, for all $t\in(0,T)$, and all $u,v\in H$.\\
(A.5) There exists a constant $\hat{K}$ such that $|G(t,u)-G(t,v)|^{2}_{L_{Q}^{V}}\leq \hat{K}||u-v||^{2}$, for all $t\in(0,T)$, and all $u,v\in V$.\\
Under the assumptions (A.1)-(A.4), we consider the following 3D stochastic Navier-Stokes equation with damping
\begin{align}
du(t)+[\mathcal{A}u(t)+B(u(t))+g(u(t))]dt&=G(t,u(t))dW(t),\\
u(0,x)&=\xi,
\end{align}
 which has a unique solution $u\in L^{2}(\Omega;C([0,T];H))\cap L^{2}(\Omega\times[0,T];V)\cap L^{\beta+1}(\Omega\times[0,T];L^{\beta+1})$ and $\nabla u\in L^{2}([0,T]\times\Omega;V)\cap L^{2}(\Omega,C([0,T];H))$ for $\beta > 3$ with any $\alpha>0$ and $\alpha \geq \frac12$ as $\beta = 3$ and we have
\begin{eqnarray}
u(t)=\xi-\int_{0}^{t}\mathcal{A}u(s)ds-\int_{0}^{t}B(u(s))ds-\int_{0}^{t}g(u(s))ds+\int_{0}^{t}G(s,u(s))dW(s).
\end{eqnarray}
Consider the small time process $u(\varepsilon t)$, we have the following equation:
\begin{eqnarray}
u^{\varepsilon}(t)=\xi-\varepsilon\int_{0}^{t}\mathcal{A}u^{\varepsilon}(s)ds-\varepsilon\int_{0}^{t}
B(u^{\varepsilon}(s))ds-\varepsilon\int_{0}^{t}g(u^{\varepsilon}(s))ds+\sqrt{\varepsilon}\int_{0}^{t}G(\varepsilon s,.)dW(s).
\end{eqnarray}
Let $\mu_{\xi}^{\varepsilon}$ be the law of $u^{\varepsilon}$ and $\nabla u^{\varepsilon}$ on $C([0,1];H)$. Define the following functional $I(g)$ by
\begin{eqnarray*}
I(g)=\inf_{h\in\Gamma_{g}}\{\frac{1}{2}\int_{0}^{1}|\dot{h}(t)|^{2}_{H_{0}}dt\},
\end{eqnarray*}
where $H_{0}=Q^{\frac{1}{2}}H$ endowed with the norm $|h|^{2}_{H_{0}}=|Q^{-\frac{1}{2}}h|^{2}$, and
\begin{eqnarray*}
\Gamma_{g}=\{h=(h_{1},h_{2},\cdots,h_{k},\cdots)\in C([0,T];H):h(\cdot)~is~absolutely~continuous~and~such~that
\end{eqnarray*}
\begin{eqnarray*}
g(t)=\xi+\int^{t}_{0}G(s,g(s))\dot{h}(s)ds,0\leq t\leq 1\}.
\end{eqnarray*}
{\bf Theorem 5.1} Assume (A.1)-(A.5) are satisfied, $\xi\in H$, then there exists a sequence $\{\xi_{n}\}_{n=1}^{\infty}\subset V$ such that $\xi_{n}\rightarrow \xi$,  for $\beta > 3$ with any $\alpha>0$ and $\alpha \geq \frac12$ as $\beta = 3$, then $\mu_{\xi}^{\varepsilon}$ satisfies a large deviation with rate function $I(\cdot)$, i.e.\\
(i) for every closed subset $F\subset C([0,1];H)$,
\begin{eqnarray*}
\lim_{\varepsilon\rightarrow 0}\sup_{\xi_{n}\rightarrow \xi}\varepsilon log\mu_{\xi_{n}}^{\varepsilon}(F)\leq -\inf_{g\in F}(I(g)),
\end{eqnarray*}
(ii) for every open subset $E\subset C([0,1];H)$,
\begin{eqnarray*}
\lim_{\varepsilon\rightarrow 0}\inf_{\xi_{n}\rightarrow \xi}\varepsilon log\mu_{\xi_{n}}^{\varepsilon}(E)\geq -\inf_{g\in E}(I(g)).
\end{eqnarray*}
Proof. This section is devoted to the proof of Theorem 5.1, which will be split into a number of lemmas. Let $v^{\varepsilon}(\cdot)$ be the solution of the stochastic equation
\begin{eqnarray}
v^{\varepsilon}(t)=\xi+\sqrt{\varepsilon}\int_{0}^{t}G(\varepsilon s,v^{\varepsilon}(s))dW(s),
\end{eqnarray}
and $\nu^{\varepsilon}$ be the law of $v^{\varepsilon}(\cdot)$ on the $C([0,1];H)$. By \cite{da2}, we obtain that $\nu^{\varepsilon}$ satisfies a large deviation principle with the rate function $I(\cdot)$. First, we will prove that the two families of probability measure $\mu^{\varepsilon}$ and $\nu^{\varepsilon}$ are exponentially equivalent, that is, for any $\delta>0$, we deduce
\begin{eqnarray}
\lim_{\varepsilon\rightarrow 0}\varepsilon log P(\sup_{0\leq t\leq 1}|u^{\varepsilon}(t)-v^{\varepsilon}(t)|^{2}>\delta)=-\infty.
\end{eqnarray}
Then Theorem 5.1 follows from (5.6) and Theorem 4.2.13 in \cite{dembo} for $\xi$. The Theorem 4.2.13 in \cite{dembo} satisfies  that if one of the two exponentially equivalent families satisfies a large deviation principle, so does the other.\\
The following result is an useful estimate of the probability that the solution of (5.4) leave an energy ball. It will play an important role in the rest of the paper.\\
{\bf Lemma 5.1}
Let $u^{\varepsilon}(t)$ be the solution of (5.4), then we have
\begin{eqnarray}
\lim_{M\rightarrow \infty}\sup_{0<\varepsilon\leq 1 }\varepsilon log P((|u^{\varepsilon}|^{H}_{V}(1))^{2}>M)=-\infty,
\end{eqnarray}
where $(|u^{\varepsilon}|^{H}_{V}(1))^{2}=\sup_{0\leq t\leq 1}|u^{\varepsilon}(t)|^{2}+2\varepsilon\int_{0}^{1}||u^{\varepsilon}(t)||^{2}dt$.\\
Proof. Applying It$\hat{o}$'s formula to $|u^{\varepsilon}(t)|^{2}$, we deduce
\begin{align*}
|u^{\varepsilon}(t)|^{2}&=|\xi|^{2}-2\varepsilon\int_{0}^{t}\langle u^{\varepsilon}(s),\mathcal{A}u^{\varepsilon}(s)\rangle ds-2\varepsilon\int_{0}^{t}\langle u^{\varepsilon}(s),B(u^{\varepsilon}(s))\rangle ds-2\varepsilon\int_{0}^{t}\langle u^{\varepsilon}(s),g(u^{\varepsilon}(s))\rangle ds\\
&+2\sqrt{\varepsilon}\int_{0}^{t}(u^{\varepsilon}(s),G(\varepsilon s,u^{\varepsilon}(s))dW(s))+\varepsilon\int_{0}^{t}|G(\varepsilon s,u^{\varepsilon}(s))|^{2}_{L_{Q}}ds,
\end{align*}
that is,
\begin{align}
&|u^{\varepsilon}(t)|^{2}+2\varepsilon\int_{0}^{t}( u^{\varepsilon}(s),\mathcal{A}u^{\varepsilon}(s)) ds+2\varepsilon\int_{0}^{t}(u^{\varepsilon}(s),g(u^{\varepsilon}(s)))ds\notag\\
&=|\xi|^{2}-2\varepsilon\int_{0}^{t}( u^{\varepsilon}(s),B(u^{\varepsilon}(s))) ds\notag\\
&+2\sqrt{\varepsilon}\int_{0}^{t}(u^{\varepsilon}(s),G(\varepsilon s,u^{\varepsilon}(s))dW(s))+\varepsilon\int_{0}^{t}|G(\varepsilon s,u^{\varepsilon}(s))|^{2}_{L_{Q}}ds \notag\\
&= |\xi|^{2}+I_{1}+I_{2}+I_{3}.
\end{align}
For the first term, we deduce
\begin{eqnarray}
I_{1}(t)=0.
\end{eqnarray}
By virtue of (A.2), we have
\begin{eqnarray}
|I_{3}(t)|\leq \varepsilon L\int_{0}^{t}(1+|u^{\varepsilon}(s)|^{2})ds.
\end{eqnarray}
Substituting (5.9) and (5.10) into (5.8),  we deduce
\begin{align*}
&|u^{\varepsilon}(t)|^{2}+2\varepsilon\int_{0}^{t}||u^{\varepsilon}(s)||^{2}ds+2\varepsilon\int_{0}^{t}|u^{\varepsilon}(s)|^{\beta+1}_{\beta+1}ds\\
&\leq(|\xi|^{2}+\varepsilon Lt)+\varepsilon L\int_{0}^{t}|u^{\varepsilon}(s)|^{2}ds\\
&+2\sqrt{\varepsilon}|\int_{0}^{t}(u^{\varepsilon}(s),G(\varepsilon s,u^{\varepsilon}(s))dW(s))|.
\end{align*}
Hence, we have the following inequality
\begin{align*}
&(|u^{\varepsilon}|^{H}_{V}(T))^{2}+2\varepsilon\int_{0}^{T}|u^{\varepsilon}(s)|^{\beta+1}_{\beta+1}ds\\
&\leq2(|\xi|^{2}+\varepsilon LT)+4\varepsilon L\int_{0}^{T}(|u^{\varepsilon}(s)|^{H}_{V}(s))^{2}ds\\
&+4\sqrt{\varepsilon}\sup_{0\leq t\leq T}|\int_{0}^{t}(u^{\varepsilon}(s),G(\varepsilon s,u^{\varepsilon}(s))dW(s))|.
\end{align*}
For $p\geq 2$, there exists a constant $c$ such that
\begin{align}
&(\mathbb{E}(|u^{\varepsilon}|^{H}_{V}(T))^{2p})^{\frac{1}{p}}\notag\\
&\leq2(|\xi|^{2}+\varepsilon LT)+c(\mathbb{E}(\int_{0}^{T}(|u^{\varepsilon}(s)|^{H}_{V}(s))^{2}ds)^{p})^{\frac{1}{p}}\notag\\
&+4\sqrt{\varepsilon}(\mathbb{E}(\sup_{0\leq t\leq T}|\int_{0}^{t}(u^{\varepsilon}(s),G(\varepsilon s,u^{\varepsilon}(s))dW(s))|^{p}))^{\frac{1}{p}}.
\end{align}
For estimating the stochastic integral term, we will use the following result from \cite{barlow, davis} that there exists a positive constant $C$ such that, for any $p\geq 2$ and for any continuous martingale $(M_{t})$ with $M_{0}=0$, one has
\begin{eqnarray}
||M^{*}_{t}||_{p}\leq Cp^{\frac{1}{2}}||\langle M\rangle^{\frac{1}{2}}_{t}||_{p},
\end{eqnarray}
where $M^{*}_{t}=\sup_{0\leq s\leq t}|M_{s}|$ and $||\cdot||_{p}$ stands for the $L^{p}-$norm.
\begin{align}
&4\sqrt{\varepsilon}(\mathbb{E}(\sup_{0\leq t\leq T}|\int_{0}^{t}(u^{\varepsilon}(s),G(\varepsilon s,u^{\varepsilon}(s))dW(s))|^{p}))^{\frac{1}{p}}\notag\\
&\leq4C\sqrt{p\varepsilon}((\mathbb{E}(\int_{0}^{T}(1+|u^{\varepsilon}(s)|^{2})^{2}ds)^{\frac{p}{2}})^{\frac{2}{p}})^{\frac{1}{2}}\notag\\
&\leq4C\sqrt{p\varepsilon}(\int_{0}^{T}1+(\mathbb{E}|u^{\varepsilon}(s)|^{2p})^{\frac{2}{p}}ds)^{\frac{1}{2}},
\end{align}
where (A.2) has been used. In other words,
\begin{eqnarray}
c(\mathbb{E}(\int_{0}^{T}(|u^{\varepsilon}|_{V}^{H}(s))^{2}ds)^{p})^{\frac{1}{p}}\leq c\int_{0}^{T}(\mathbb{E}(|u^{\varepsilon}|^{H}_{V}(s))^{2p})^{\frac{1}{p}}ds.
\end{eqnarray}
Combining (5.11), (5.13) and (5.14), we obtain
\begin{align}
&(\mathbb{E}(|u^{\varepsilon}|^{H}_{V}(T))^{2p})^{\frac{2}{p}}\notag\\
&\leq 8(|\xi|^{2}+\varepsilon LT)^{2}+2c^{2}T \int_{0}^{T}(\mathbb{E}(|u^{\varepsilon}(s)|^{H}_{V}(s))^{2p})^{\frac{2}{p}}ds\notag\\
&+ 32C^{2}p\varepsilon T+32C^{2}p\varepsilon\int_{0}^{T}(\mathbb{E}(|u^{\varepsilon}(s)|^{H}_{V}(s))^{2p})^{\frac{2}{p}}ds.
\end{align}
Applying the Gronwall's inequality, we get
\begin{eqnarray}
(\mathbb{E}(|u^{\varepsilon}|^{H}_{V}(1))^{2p})^{\frac{2}{p}}\leq C[8(|\xi|^{2}+\varepsilon L)^{2}+32C^{2}p\varepsilon]exp(2c^{2}+32C^{2}p\varepsilon).
\end{eqnarray}
Since $P((|u^{\varepsilon}|^{H}_{V}(1))^{2}> M)\leq M^{-p}E(|u^{\varepsilon}|^{H}_{V}(1))^{2p}$. Let $p=\frac{2}{\varepsilon}$ in (5.16), we obtain
\begin{align*}
\varepsilon logP((|u^{\varepsilon}|^{H}_{V}(1))^{2}> M)&\leq -logM+log(E(|u^{\varepsilon}|^{H}_{V}(1))^{2p})^{\frac{1}{p}}\notag\\
&\leq -logM+log\sqrt{C[8(|\xi|^{2}+\varepsilon L)^{2}+32C^{2}]}+c^{2}+32C^{2}.
\end{align*}
Hence, we get the following inequality
\begin{eqnarray*}
\sup_{0<\varepsilon\leq 1}\varepsilon logP((|u^{\varepsilon}|^{H}_{V}(1))^{2}> M)\leq -logM+log\sqrt{C[8(|\xi|^{2}+ L)^{2}+32C^{2}]}+c^{2}+32C^{2}.
\end{eqnarray*}
Letting $M\rightarrow \infty$ on both side of the above inequality, this completes the proof of Lemma 5.1.\\
For $\xi\in H$, since $V$ is dense in $H$, then there exists a sequence $\{\xi_{n}\}_{n=1}^{\infty}\subset V$ such that
\begin{eqnarray*}
\lim_{n\rightarrow +\infty}|\xi_{n}-\xi|=0.
\end{eqnarray*}
Let $u_{n}^{\varepsilon}(\cdot)$ be the solution of (5.4) with the initial value $\xi_{n}$. By the proof of Lemma 5.1, it is easy to see that we get the following equality
\begin{eqnarray}
\lim_{M\rightarrow +\infty}\sup_{n}\sup_{0<\varepsilon\leq 1}\varepsilon logP((|u^{\varepsilon}_{n}|^{H}_{V}(1))^{2}> M)=-\infty.
\end{eqnarray}
Let $v_{n}^{\varepsilon}(\cdot)$ be the solution of (5.5) with the initial value $\xi_{n}$, and we have the following result whose proof is similar to the Lemma 5.1.\\
{\bf Lemma 5.2} For any fixed $n\in Z^{+}$,
\begin{eqnarray*}
\lim_{M\rightarrow +\infty}\sup_{0<\varepsilon\leq 1}\varepsilon logP(\sup_{0\leq t\leq 1}||v^{\varepsilon}_{n}(t)||^{2}>M)=-\infty,\\
\lim_{M\rightarrow +\infty}\sup_{0<\varepsilon\leq 1}\varepsilon logP(\sup_{0\leq t\leq 1}| v^{\varepsilon}_{n}(t)|^{\beta+1}_{\beta+1}>M)=-\infty.
\end{eqnarray*}
Proof. Applying It$\hat{o}$'s formula to $||v_{n}^{\varepsilon}(t)||^{2}$, we deduce
\begin{align*}
||v_{n}^{\varepsilon}(t)||^{2}&=||\xi_{n}||^{2}
+2\sqrt{\varepsilon}\int_{0}^{t}((v_{n}^{\varepsilon}(s),G(\varepsilon s,v_{n}^{\varepsilon}(s))dW(s)))
+\varepsilon\int^{t}_{0}|G(\varepsilon s,v_{n}^{\varepsilon}(s))|^{2}_{L^{V}_{Q}}ds.
\end{align*}
By (A.3) and inequality (5.12), we obtain
\begin{align}
(\mathbb{E}[\sup_{0\leq t\leq r}||v_{n}^{\varepsilon}(t)||^{2p}])^{\frac{2}{p}}&\leq 2||\xi_{n}||^{4}+8C\varepsilon p(\mathbb{E}[\int_{0}^{r}||v_{n}^{\varepsilon}(s)||^{2}|G(\varepsilon s,v_{n}^{\varepsilon}(s))|^{2}_{L^{V}_{Q}}ds]^{\frac{p}{2}})^{\frac{2}{p}}\notag\\
&+ 4\varepsilon^{2}\hat{L}^{2}r(r+\int_{0}^{r}(\mathbb{E}[\sup_{0\leq l\leq s}||v_{n}^{\varepsilon}(l)||^{2p}])^{\frac{2}{p}}ds)\notag\\
&\leq 2||\xi_{n}||^{4}+16C\varepsilon p\hat{L}(r+\int_{0}^{r}(\mathbb{E}[\sup_{0\leq l\leq s}||v_{n}^{\varepsilon}(l)||^{2p}])^{\frac{2}{p}}ds)\notag\\
&+ 4\varepsilon^{2}\hat{L}^{2}r(r+\int_{0}^{r}(\mathbb{E}[\sup_{0\leq l\leq s}||v_{n}^{\varepsilon}(l)||^{2p}])^{\frac{2}{p}}ds).
\end{align}
By Gronwall's inequality,
\begin{eqnarray}
(\mathbb{E}[\sup_{0\leq t\leq1}||v_{n}^{\varepsilon}(t)||^{2p}])^{\frac{2}{p}}\leq (2||\xi_{n}||^{4}+16C\varepsilon p\hat{L}+4\varepsilon^{2}\hat{L}^{2})e^{16C\varepsilon p\hat{L}+4\varepsilon^{2}\hat{L}^{2}}.
\end{eqnarray}
According to the proof of Lemma 5.1,  we have
\begin{eqnarray*}
\lim_{M\rightarrow +\infty}\sup_{0<\varepsilon\leq 1}\varepsilon logP(\sup_{0\leq t\leq 1}||v^{\varepsilon}_{n}(t)||^{2}>M)=-\infty.
\end{eqnarray*}
Similarly, we can easy to get the following equality
\begin{eqnarray*}
\lim_{M\rightarrow +\infty}\sup_{0<\varepsilon\leq 1}\varepsilon logP(\sup_{0\leq t\leq 1}| v^{\varepsilon}_{n}(t)|^{\beta+1}_{\beta+1}>M)=-\infty.
\end{eqnarray*}
{\bf Lemma 5.3} For any $\delta>0$ and $\beta > 3$ with any $\alpha>0$ and $\alpha \geq \frac12$ as $\beta = 3$,
\begin{eqnarray}
\lim_{n\rightarrow\infty}\sup_{0<\varepsilon\leq 1}\varepsilon logP(\sup_{0\leq t\leq 1}|u^{\varepsilon}(t)-u_{n}^{\varepsilon}(t)|^{2}>\delta)=-\infty.
\end{eqnarray}
Proof. For any $M>0$, define a stopping time
\begin{eqnarray*}
\tau_{\varepsilon,M}=inf\{t:\varepsilon\int_{0}^{t}|| u^{\varepsilon}(r)||^{2}dr>M,~or~|u^{\varepsilon}(t)|^{2}>M,or~||u^{\varepsilon}(t)||^{2}>M\}.
\end{eqnarray*}
We deduce
\begin{align}
&P(\sup_{0\leq t\leq 1}|u^{\varepsilon}(t)-u_{n}^{\varepsilon}(t)|^{2}>\delta,(|u^{\varepsilon}|^{H}_{V}(1))^{2}\leq M,||u^{\varepsilon}(t)||^{2}\leq M)\notag\\
&\leq P(\sup_{0\leq t\leq 1}|u^{\varepsilon}(t)-u_{n}^{\varepsilon}(t)|^{2}>\delta,\tau_{\varepsilon,M}\geq1)\notag\\
&\leq P(\sup_{0\leq t\leq 1}|u^{\varepsilon}(t)-u_{n}^{\varepsilon}(t)|^{2}>\delta).
\end{align}
Let $k$ be a positive constant. Applying It$\hat{o}$'s formula to $e^{-k\varepsilon\int_{0}^{t\wedge\tau_{\varepsilon,M}}||
u^{\varepsilon}(s)||^{4}ds}|u^{\varepsilon}(t\wedge\tau_{\varepsilon,M})-u_{n}^{\varepsilon}
(t\wedge\tau_{\varepsilon,M})|^{2}$, we deduce
\begin{align}
&e^{-k\varepsilon\int_{0}^{t\wedge\tau_{\varepsilon,M}}|| u^{\varepsilon}||^{4}ds}|u^{\varepsilon}(t\wedge\tau_{\varepsilon,M})-u_{n}^{\varepsilon}
(t\wedge\tau_{\varepsilon,M})|^{2}\notag\\
&+2\varepsilon\int_{0}^{t\wedge\tau_{\varepsilon,M}}e^{-k\varepsilon\int_{0}^{s}||u^{\varepsilon}||^{4}dr}
||u^{\varepsilon}(s)-u^{\varepsilon}_{n}(s)||^{2}ds\notag\\
&=|\xi-\xi_{n}|^{2}-k\varepsilon\int_{0}^{t\wedge\tau_{\varepsilon,M}}e^{-k\varepsilon\int_{0}^{s}||u^{\varepsilon}(r)||^{4}dr}
||u^{\varepsilon}(s)||^{4}|u^{\varepsilon}(s)-u_{n}^{\varepsilon}(s)|^{2}ds\notag\\
&- 2\varepsilon\int_{0}^{t\wedge\tau_{\varepsilon,M}}e^{-k\varepsilon\int_{0}^{s}||u^{\varepsilon}(r)||^{4}dr}(b(u^{\varepsilon}(s),u^{\varepsilon}(s),
u^{\varepsilon}(s)-u^{\varepsilon}_{n}(s))\notag\\
&-b(u^{\varepsilon}_{n}(s),u^{\varepsilon}_{n}(s),
u^{\varepsilon}(s)-u^{\varepsilon}_{n}(s)))ds\notag\\
&-2\varepsilon\int_{0}^{t\wedge\tau_{\varepsilon,M}}e^{-k\varepsilon\int_{0}^{s}|| u^{\varepsilon}(r)||^{4}dr}(g(u^{\varepsilon}(s))-g(u_{n}^{\varepsilon}(s)),u^{\varepsilon}(s)-u_{n}^{\varepsilon}(s))\notag\\
&+2\sqrt{\varepsilon}\int_{0}^{t\wedge\tau_{\varepsilon,M}}e^{-k\varepsilon\int_{0}^{s}|| u^{\varepsilon}(r)||^{4}dr}(u^{\varepsilon}(s)-u_{n}^{\varepsilon}(s),(G(\varepsilon s,u^{\varepsilon}(s))-G(\varepsilon s,u^{\varepsilon}_{n}(s)))dW(s))\notag\\
&+\varepsilon\int_{0}^{t\wedge\tau_{\varepsilon,M}}e^{-k\varepsilon\int_{0}^{s}||u^{\varepsilon}(r)||^{4}dr}|G(\varepsilon s,u^{\varepsilon}(s))-G(\varepsilon s,u^{\varepsilon}_{n}(s))|^{2}_{L_{Q}}ds.
\end{align}
First, we obtain
\begin{eqnarray*}
b(u_{n}^{\varepsilon}(t),u_{n}^{\varepsilon}(t),u^{\varepsilon}(t)-u_{n}^{\varepsilon}(t))=b(u_{n}^{\varepsilon}(t),u^{\varepsilon}(t),u^{\varepsilon}(t)-
u_{n}^{\varepsilon}(t)),
\end{eqnarray*}
and
\begin{align*}
&|b(u^{\varepsilon}(t),u^{\varepsilon}(t),u^{\varepsilon}(t)-u_{n}^{\varepsilon}(t))-b(u_{n}^{\varepsilon}(t),u^{\varepsilon}_{n}(t),u^{\varepsilon}(t)-
u_{n}^{\varepsilon}(t))|\\
&=|b(u^{\varepsilon}(t)-u_{n}^{\varepsilon}(t),u^{\varepsilon}(t),u^{\varepsilon}(t)-u_{n}^{\varepsilon}(t))|\\
&\leq c|u^{\varepsilon}(t)-u_{n}^{\varepsilon}(t)|^{\frac{1}{2}}\cdot||u^{\varepsilon}(t)||\cdot||u^{\varepsilon}(t)-u_{n}^{\varepsilon}(t)||^{\frac{3}{2}}.
\end{align*}
Combining (5.22), we deduce
\begin{align}
&e^{-k\varepsilon\int_{0}^{t\wedge\tau_{\varepsilon,M}}||u^{\varepsilon}(s)||^{4}ds}|u^{\varepsilon}(t\wedge\tau_{\varepsilon,M})-u_{n}^{\varepsilon}
(t\wedge\tau_{\varepsilon,M})|^{2}\notag\\
&+2\varepsilon\int_{0}^{t\wedge\tau_{\varepsilon,M}}e^{-k\varepsilon\int_{0}^{s}||u^{\varepsilon}(r)||^{4}dr}
||u^{\varepsilon}(s)-u^{\varepsilon}_{n}(s)||^{2}ds\notag\\
&\leq|\xi-\xi_{n}|^{2}-k\varepsilon\int_{0}^{t\wedge\tau_{\varepsilon,M}}e^{-k\varepsilon\int_{0}^{s}||u^{\varepsilon}(r)||^{4}dr}
|| u^{\varepsilon}(s)||^{4}|u^{\varepsilon}(s)-u_{n}^{\varepsilon}(s)|^{2}ds\notag\\
&+ 2c\varepsilon\int_{0}^{t\wedge\tau_{\varepsilon,M}}e^{-k\varepsilon\int_{0}^{s}||u^{\varepsilon}(r)||^{4}dr}|u^{\varepsilon}(t)-u_{n}^{\varepsilon}(t)
|^{\frac{1}{2}}\cdot||u^{\varepsilon}(t)||\cdot||u^{\varepsilon}(t)-u_{n}^{\varepsilon}(t)||^{\frac{3}{2}}ds\notag\\
&+2\sqrt{\varepsilon}|\int_{0}^{t\wedge\tau_{\varepsilon,M}}e^{-k\varepsilon\int_{0}^{s}|| u^{\varepsilon}(r)||^{4}dr}(u^{\varepsilon}(s)-u_{n}^{\varepsilon}(s),
(G(\varepsilon s,u^{\varepsilon}(s))-G(\varepsilon s,u^{\varepsilon}_{n}(s)))dW(s))|\notag\\
&+K\varepsilon\int_{0}^{t\wedge\tau_{\varepsilon,M}}e^{-k\varepsilon\int_{0}^{s}||u^{\varepsilon}(r)||^{4}dr}|u^{\varepsilon}(s)-u^{\varepsilon}_{n}(s)|^{2}ds.
\notag\\
&\leq|\xi-\xi_{n}|^{2}-k\varepsilon\int_{0}^{t\wedge\tau_{\varepsilon,M}}e^{-k\varepsilon\int_{0}^{s}||u^{\varepsilon}(r)||^{4}dr}
||u^{\varepsilon}(s)||^{4}|u^{\varepsilon}(s)-u_{n}^{\varepsilon}(s)|^{2}ds\notag\\
&+ C\varepsilon\int_{0}^{t\wedge\tau_{\varepsilon,M}}e^{-k\varepsilon\int_{0}^{s}||u^{\varepsilon}(r)||^{4}dr}|u^{\varepsilon}(s)-u_{n}^{\varepsilon}(s)|^{2}
\cdot||u_{n}^{\varepsilon}(s)||^{4}ds\notag\\
&+\varepsilon\int_{0}^{t\wedge\tau_{\varepsilon,M}}e^{-k\varepsilon\int_{0}^{s}||u^{\varepsilon}(r)||^{4}dr}||u^{\varepsilon}(s)-u_{n}^{\varepsilon}(s)
||^{2}ds\notag\\
&+2\sqrt{\varepsilon}|\int_{0}^{t\wedge\tau_{\varepsilon,M}}e^{-k\varepsilon\int_{0}^{s}|| u^{\varepsilon}(r)||^{4}dr}(u^{\varepsilon}(s)-u_{n}^{\varepsilon}(s),
(G(\varepsilon s,u^{\varepsilon}(s))-G(\varepsilon s,u^{\varepsilon}_{n}(s)))dW(s))|\notag\\
&+K\varepsilon\int_{0}^{t\wedge\tau_{\varepsilon,M}}e^{-k\varepsilon\int_{0}^{s}||u^{\varepsilon}(r)||^{4}dr}|u^{\varepsilon}(s)-u^{\varepsilon}_{n}(s)|^{2}ds.
\end{align}
Choosing $k>C$ and using (5.12), we deduce
\begin{align}
&(\mathbb{E}[\sup_{0\leq s\leq t\wedge\tau_{\varepsilon,M}}(e^{-k\varepsilon\int_{0}^{s}||u^{\varepsilon}(r)||^{4}dr}|u^{\varepsilon}(s)
-u_{n}^{\varepsilon}(s)|^{2})]^{p})^{\frac{2}{p}}\notag\\
&\leq 2|\xi-\xi_{n}|^{4}+2\varepsilon^{2}K^{2}\int_{0}^{t}(\mathbb{E}[(\sup_{0\leq r\leq s\wedge\tau_{\varepsilon,M}}e^{-k\varepsilon\int_{0}^{r}|| u^{\varepsilon}(l)||^{4}dl}|u^{\varepsilon}(r)
-u_{n}^{\varepsilon}(r)|^{2})^{p}])^{\frac{2}{p}}ds\notag\\
&+8C\varepsilon pK^{2}\int_{0}^{t}(\mathbb{E}[(\sup_{0\leq r\leq s\wedge\tau_{\varepsilon,M}}e^{-2k\varepsilon\int_{0}^{r}|| u^{\varepsilon}(l)||^{4}dl}|u^{\varepsilon}(r)
-u_{n}^{\varepsilon}(r)|^{4})^{\frac{p}{2}}])^{\frac{2}{p}}ds\notag\\
&\leq 2|\xi-\xi_{n}|^{4}+2\varepsilon^{2}K^{2}\int_{0}^{t}(\mathbb{E}[(\sup_{0\leq r\leq s\wedge\tau_{\varepsilon,M}}e^{-k\varepsilon\int_{0}^{r}|| u^{\varepsilon}(l)||^{4}dl}|u^{\varepsilon}(r)
-u_{n}^{\varepsilon}(r)|^{2})^{p}])^{\frac{2}{p}}ds\notag\\
&+8C\varepsilon pK^{2}\int_{0}^{t}(\mathbb{E}[(\sup_{0\leq r\leq s\wedge\tau_{\varepsilon,M}}e^{-k\varepsilon\int_{0}^{r}|| u^{\varepsilon}(l)||^{4}dl}|u^{\varepsilon}(r)
-u_{n}^{\varepsilon}(r)|^{2})^{p}])^{\frac{2}{p}}ds.
\end{align}
Applying Gronwall's inequality, we have the following inequality
\begin{align}
&(\mathbb{E}[\sup_{0\leq t\leq 1\wedge\tau_{\varepsilon,M}}(e^{-k\varepsilon\int_{0}^{t}||u^{\varepsilon}(s)||^{4}ds}|u^{\varepsilon}(t)
-u_{n}^{\varepsilon}(t)|^{2})^{p}])^{\frac{2}{p}}\notag\\
&\leq 2|\xi-\xi_{n}|^{4}e^{2\varepsilon^{2}K^{2}+8C\varepsilon pK^{2}}.
\end{align}
Hence, we have
\begin{align}
&(\mathbb{E}[\sup_{0\leq t\leq 1\wedge\tau_{\varepsilon,M}}(|u^{\varepsilon}(t)-u_{n}^{\varepsilon}(t)|^{2})^{p}])^{\frac{2}{p}}\notag\\
&\leq (\mathbb{E}[\sup_{0\leq t\leq 1\wedge\tau_{\varepsilon,M}}(e^{-k\varepsilon\int_{0}^{t}||u^{\varepsilon}(s)||^{4}ds}|u^{\varepsilon}(t)
-u_{n}^{\varepsilon}(t)|^{2})^{p}e^{kp\varepsilon\int_{0}^{1\wedge\tau_{\varepsilon,M}}||u^{\varepsilon}(s)||^{4}ds}])^{\frac{2}{p}}\notag\\
&\leq e^{2kM}(\mathbb{E}[\sup_{0\leq t\leq 1\wedge\tau_{\varepsilon,M}}(e^{-k\varepsilon\int_{0}^{t}||u^{\varepsilon}(s)||^{4}ds}|u^{\varepsilon}(t)
-u_{n}^{\varepsilon}(t)|^{2})^{p}])^{\frac{2}{p}}\notag\\
&\leq 2e^{2kM}|\xi-\xi_{n}|^{4}e^{2\varepsilon^{2}K^{2}+8C\varepsilon pK^{2}}.
\end{align}
Fix a positive constant $M$ and $p=\frac{2}{\varepsilon}$, we have the following inequality
\begin{align}
&\sup_{0<\varepsilon\leq 1}\varepsilon log P(\sup_{0\leq t\leq 1\wedge\tau_{\varepsilon,M}}(|u^{\varepsilon}(t)-u_{n}^{\varepsilon}(t)|^{2}>\delta)\notag\\
&\leq \sup_{0<\varepsilon\leq 1}\varepsilon log\frac{\mathbb{E}[\sup_{0\leq t\leq 1\wedge\tau_{\varepsilon,M}}|u^{\varepsilon}(t)-u_{n}^{\varepsilon}(t)|^{2p}]}{\delta^{p}}\notag\\
&\leq 2k M+2K^{2}+16CK^{2}-2log\delta+log2|\xi-\xi_{n}|^{4}\notag\\
&\rightarrow -\infty,~~as~n\rightarrow +\infty.
\end{align}
For any given $R>0$, by Lemma 5.1, there exists a positive constant $M$ such that for any $\varepsilon\in(0,1]$, we get the following inequalities
\begin{eqnarray}
P((|u^{\varepsilon}|^{H}_{V}(1))^{2}>M)\leq e^{-\frac{R}{\varepsilon}}.
\end{eqnarray}
For any $M$, (5.22) and (5.27) implies that there exists a positive integer $N$, such that for any $n\geq N$,
\begin{eqnarray}
\sup_{0<\varepsilon\leq 1}\varepsilon log P(\sup_{0\leq t\leq 1}|u^{\varepsilon}(t)-u_{n}^{\varepsilon}(t)|^{2}>\delta,(|u^{\varepsilon}|^{H}_{V}(1))^{2}\leq M, ||u^{\varepsilon}(t)||^{2}\leq M)
\leq-R.
\end{eqnarray}
Combining (5.28) and (5.29) together, we can get a positive number $N$, such that for any $n\geq N$, $\varepsilon\in(0,1],$
\begin{eqnarray}
P(\sup_{0\leq t\leq 1}|u^{\varepsilon}(t)-u_{n}^{\varepsilon}(t)|^{2}>\delta)\leq 2e^{-\frac{R}{\varepsilon}}.
\end{eqnarray}
Since $R$ is arbitrary, this completes the proof of Lemma 5.3.

The following lemma can be proved similarly as Lemma 5.3.\\
{\bf Lemma 5.4} For any $\delta>0$,
\begin{eqnarray}
\lim_{n\rightarrow\infty}\sup_{0<\varepsilon\leq 1}\varepsilon logP(\sup_{0\leq t\leq 1}|v^{\varepsilon}(t)-v_{n}^{\varepsilon}(t)|^{2}>\delta)=-\infty.
\end{eqnarray}
The following lemma have that for any fixed integer $n$, the two families $\{u_{n}^{\varepsilon},\varepsilon>0\}$ and $\{v_{n}^{\varepsilon},\varepsilon>0\}$ are exponentially equivalent.\\
{\bf Lemma 5.5} For any $\delta>0$ and every positive integer $n$ and $\beta > 3$ with any $\alpha>0$ and $\alpha \geq \frac12$ as $\beta = 3$,
\begin{eqnarray}
\lim_{\varepsilon\rightarrow 0}\varepsilon logP(\sup_{0\leq t\leq 1}|u_{n}^{\varepsilon}(t)-v_{n}^{\varepsilon}(t)|^{2}>\delta)=-\infty.
\end{eqnarray}
Proof. For any $M>0$, define stopping times
\begin{align*}
\tau^{1,n}_{\varepsilon,M}&=inf\{t\geq 0; \varepsilon\int_{0}^{t}||u_{n}^{\varepsilon}(s)||^{2}ds>M,~or~|u_{n}^{\varepsilon}(t)|^{2}>M,~or~||u_{n}^{\varepsilon}(s)||^{2}>M\},\\
\tau^{2,n}_{\varepsilon,M}&=inf\{t\geq 0; ||v_{n}^{\varepsilon}(t)||^{2}>M,| v_{n}^{\varepsilon}(t)|^{\beta+1}_{\beta+1}>M\}.
\end{align*}
We deduce
\begin{align}
&P(\sup_{0\leq t\leq1}|u_{n}^{\varepsilon}(t)-v_{n}^{\varepsilon}(t)|^{2}>\delta,(|u_{n}^{\varepsilon}|_{V}^{H}(1))^{2}\leq M,\sup_{0\leq t\leq1}||u_{n}^{\varepsilon}(t)||^{2}\leq M,\notag\\
&\sup_{0\leq t\leq1}||v_{n}^{\varepsilon}(t)||^{2}\leq M,\sup_{0\leq t\leq1}| v_{n}^{\varepsilon}(t)|^{\beta+1}_{\beta+1}\leq M)\notag\\
&\leq P(\sup_{0\leq t\leq1}|u_{n}^{\varepsilon}(t)-v_{n}^{\varepsilon}(t)|^{2}>\delta,1\leq\tau^{1,n}_{\varepsilon,M}\wedge\tau^{2,n}_{\varepsilon,M})\notag\\
&\leq P(\sup_{0\leq t\leq1\wedge\tau^{1,n}_{\varepsilon,M}\wedge\tau^{2,n}_{\varepsilon,M}}|u_{n}^{\varepsilon}(t)-v_{n}^{\varepsilon}(t)|^{2}>\delta).
\end{align}
Let $\tau^{n}_{\varepsilon,M}=\tau^{1,n}_{\varepsilon,M}\wedge\tau^{2,n}_{\varepsilon,M}$. By It$\hat{o}$'s formula, we deduce
\begin{align}
&|u_{n}^{\varepsilon}(t\wedge\tau^{n}_{\varepsilon,M})-v_{n}^{\varepsilon}
(t\wedge\tau_{\varepsilon,M}^{n})|^{2}
+2\varepsilon\int_{0}^{t\wedge\tau_{\varepsilon,M}^{n}}
||u_{n}^{\varepsilon}(s)-v^{\varepsilon}_{n}(s)||^{2}ds\notag\\
&=-2\varepsilon\int_{0}^{t\wedge\tau_{\varepsilon,M}^{n}}
( u^{\varepsilon}_{n}(s)-v^{\varepsilon}_{n}(s),\mathcal{A}v^{\varepsilon}_{n}(s))ds-2\varepsilon\int_{0}^{t\wedge\tau_{\varepsilon,M}^{n}}
(u^{\varepsilon}_{n}(s)-v^{\varepsilon}_{n}(s),g(u^{\varepsilon}_{n}(s))) ds\notag\\
&- 2\varepsilon\int_{0}^{t\wedge\tau_{\varepsilon,M}^{n}}b(u^{\varepsilon}_{n}(s),u^{\varepsilon}_{n}(s)),
u^{\varepsilon}_{n}(s)-v^{\varepsilon}_{n}(s))ds\notag\\
&+2\sqrt{\varepsilon}\int_{0}^{t\wedge\tau_{\varepsilon,M}^{n}}(u_{n}^{\varepsilon}(s)-v_{n}^{\varepsilon}(s),
(G(\varepsilon s,u_{n}^{\varepsilon}(s))-G(\varepsilon s,v^{\varepsilon}_{n}(s)))dW(s))\notag\\
&+\varepsilon\int_{0}^{t\wedge\tau_{\varepsilon,M}^{n}}|G(\varepsilon s,u^{\varepsilon}_{n}(s))-G(\varepsilon s,v^{\varepsilon}_{n}(s))|^{2}_{L_{Q}}ds\notag\\
&=-2\varepsilon\int_{0}^{t\wedge\tau_{\varepsilon,M}^{n}}
( u^{\varepsilon}_{n}(s)-v^{\varepsilon}_{n}(s),\mathcal{A}v^{\varepsilon}_{n}(s))ds-2\varepsilon\int_{0}^{t\wedge\tau_{\varepsilon,M}^{n}}
(u^{\varepsilon}_{n}(s)-v^{\varepsilon}_{n}(s),g(u^{\varepsilon}_{n}(s))) ds\notag\\
&- 2\varepsilon\int_{0}^{t\wedge\tau_{\varepsilon,M}^{n}}b(u^{\varepsilon}_{n}(s)-v^{\varepsilon}_{n}(s),u^{\varepsilon}_{n}(s),
u^{\varepsilon}_{n}(s)-v^{\varepsilon}_{n}(s))ds\notag\\
&+ 2\varepsilon\int_{0}^{t\wedge\tau_{\varepsilon,M}^{n}}b(v^{\varepsilon}_{n}(s),u^{\varepsilon}_{n}(s)-v^{\varepsilon}_{n}(s),
v^{\varepsilon}_{n}(s))ds\notag\\
&+2\sqrt{\varepsilon}\int_{0}^{t\wedge\tau_{\varepsilon,M}^{n}}(u_{n}^{\varepsilon}(s)-v_{n}^{\varepsilon}(s),
(G(\varepsilon s,u_{n}^{\varepsilon}(s))-G(\varepsilon s,v^{\varepsilon}_{n}(s)))dW(s))\notag\\
&+\varepsilon\int_{0}^{t\wedge\tau_{\varepsilon,M}^{n}}|G(\varepsilon s,u^{\varepsilon}_{n}(s))-G(\varepsilon s,v^{\varepsilon}_{n}(s))|^{2}_{L_{Q}}ds\notag\\
&=I_{1}+I_{2}+I_{3}+I_{4}+I_{5}+I_{6}.
\end{align}
For the first term, we have
\begin{align}
I_{1}&\leq2\varepsilon\int_{0}^{t\wedge\tau_{\varepsilon,M}^{n}}
|( u^{\varepsilon}_{n}(s)-v^{\varepsilon}_{n}(s),\mathcal{A}v^{\varepsilon}_{n}(s))|ds\notag\\
&\leq 2\varepsilon\int_{0}^{t\wedge\tau_{\varepsilon,M}^{n}}||u^{\varepsilon}_{n}(s)-v^{\varepsilon}_{n}(s)||\cdot||v^{\varepsilon}_{n}(s)||ds.
\end{align}
For the second term, using Young inequality, we deduce
\begin{align}
I_{2}&=-2\varepsilon\int_{0}^{t\wedge\tau_{\varepsilon,M}^{n}}
(u^{\varepsilon}_{n}(s)-v^{\varepsilon}_{n}(s),g(u^{\varepsilon}_{n}(s))) ds\notag\\
&\leq -2\varepsilon\int_{0}^{t\wedge\tau_{\varepsilon,M}^{n}}|u^{\varepsilon}_{n}(s)|^{\beta+1}_{\beta+1}ds+2\varepsilon\int_{0}^{t\wedge\tau_{\varepsilon,M}^{n}}
\int_{D}|u^{\varepsilon}_{n}(s)|^{\beta}|v^{\varepsilon}_{n}(s)|dxds\notag\\
&\leq -2\varepsilon\int_{0}^{t\wedge\tau_{\varepsilon,M}^{n}}|u^{\varepsilon}_{n}(s)|^{\beta+1}_{\beta+1}ds+
2\varepsilon\int_{0}^{t\wedge\tau_{\varepsilon,M}^{n}}(\xi_{1}|u^{\varepsilon}_{n}(s)|^{\beta+1}_{\beta+1}+C(\xi_{1})|v^{\varepsilon}_{n}(s)|^{\beta+1}_{\beta+1}
)ds\notag\\
&\leq -2\varepsilon(1-\xi_{1})\int_{0}^{t\wedge\tau_{\varepsilon,M}^{n}}|u^{\varepsilon}_{n}(s)|^{\beta+1}_{\beta+1}ds+
2\varepsilon C\int_{0}^{t\wedge\tau_{\varepsilon,M}^{n}}|v^{\varepsilon}_{n}(s)|^{\beta+1}_{\beta+1}ds.
\end{align}
Next, we have
\begin{align}
I_{3}&\leq2\varepsilon\int_{0}^{t\wedge\tau_{\varepsilon,M}^{n}}|b(u^{\varepsilon}_{n}(s)-v^{\varepsilon}_{n}(s),u^{\varepsilon}_{n}(s),
u^{\varepsilon}_{n}(s)-v^{\varepsilon}_{n}(s))|ds\notag\\
&\leq c\varepsilon\int_{0}^{t\wedge\tau_{\varepsilon,M}^{n}}|u^{\varepsilon}_{n}(s)-v_{n}^{\varepsilon}(s)|^{\frac{1}{2}}\cdot||u^{\varepsilon}_{n}(s)||\cdot||u^{\varepsilon}(s)
-v_{n}^{\varepsilon}(s)||^{\frac{3}{2}}ds.
\end{align}
Similarly, we get
\begin{align}
I_{4}&\leq2\varepsilon\int_{0}^{t\wedge\tau_{\varepsilon,M}^{n}}|b(v^{\varepsilon}_{n}(s),u^{\varepsilon}_{n}(s)-v^{\varepsilon}_{n}(s),
v^{\varepsilon}_{n}(s))|ds\notag\\
&\leq c\varepsilon\int_{0}^{t\wedge\tau_{\varepsilon,M}^{n}}
|v_{n}^{\varepsilon}(s)|^{\frac{1}{2}}\cdot||u^{\varepsilon}_{n}(s)-v^{\varepsilon}_{n}(s)||\cdot||v_{n}^{\varepsilon}(s)||^{\frac{3}{2}}ds.
\end{align}
By (A.4), we have
\begin{align}
I_{6}&=\varepsilon\int_{0}^{t\wedge\tau_{\varepsilon,M}^{n}}|G(\varepsilon s,u^{\varepsilon}_{n}(s))-G(\varepsilon s,v^{\varepsilon}_{n}(s))|^{2}_{L_{Q}}ds\notag\\
&\leq \varepsilon K\int_{0}^{t\wedge\tau_{\varepsilon,M}^{n}}|u^{\varepsilon}_{n}(s)-v^{\varepsilon}_{n}(s)|^{2}ds.
\end{align}
Using (5.35)-(5.39) in (5.34) and choosing $\xi_{1}$ small enough, then there exists a constant $c$ such that
\begin{align}
&|u_{n}^{\varepsilon}(t\wedge\tau^{n}_{\varepsilon,M})-v_{n}^{\varepsilon}
(t\wedge\tau_{\varepsilon,M}^{n})|^{2}
+2\varepsilon\int_{0}^{t\wedge\tau_{\varepsilon,M}^{n}}
||u_{n}^{\varepsilon}(s)-v^{\varepsilon}_{n}(s)||^{2}ds\notag\\
&\leq\varepsilon\int_{0}^{t\wedge\tau_{\varepsilon,M}^{n}}
||u_{n}^{\varepsilon}(s)-v^{\varepsilon}_{n}(s)||^{2}ds+4\varepsilon\int_{0}^{t\wedge\tau^{n}_{\varepsilon,M}}||v^{\varepsilon}_{n}(s)||^{2}ds\notag\\
&+c\varepsilon\int_{0}^{t\wedge\tau_{\varepsilon,M}^{n}}|u^{\varepsilon}_{n}(s)-v_{n}^{\varepsilon}(s)|^{2}\cdot||u^{\varepsilon}_{n}(s)||^{4}ds \notag\\
&+2\varepsilon C\int_{0}^{t\wedge\tau_{\varepsilon,M}^{n}}|v^{\varepsilon}_{n}(s)|^{\beta+1}_{\beta+1}ds\notag\\
&+ c\varepsilon\int_{0}^{t\wedge\tau_{\varepsilon,M}^{n}}
(|v_{n}^{\varepsilon}(s)|^{2}+||v_{n}^{\varepsilon}(s)||^{6})ds+
\varepsilon K\int_{0}^{t\wedge\tau_{\varepsilon,M}^{n}}|u_{n}^{\varepsilon}(s)-v^{\varepsilon}_{n}(s)|^{2}ds\notag\\
&+2\sqrt{\varepsilon}|\int_{0}^{t\wedge\tau_{\varepsilon,M}^{n}}(u^{\varepsilon}_{n}(s)-v_{n}^{\varepsilon}(s),
(G(\varepsilon s,u^{\varepsilon}_{n}(s))-G(\varepsilon s,v^{\varepsilon}_{n}(s)))dW(s))|.
\end{align}
Applying Gronwall's inequality, we get
\begin{align}
&|u_{n}^{\varepsilon}(t\wedge\tau^{n}_{\varepsilon,M})-v_{n}^{\varepsilon}
(t\wedge\tau_{\varepsilon,M}^{n})|^{2}\notag\\
&\leq(4\varepsilon\int_{0}^{t\wedge\tau^{n}_{\varepsilon,M}}||v^{\varepsilon}_{n}(s)||^{2}ds+ c\varepsilon\int_{0}^{t\wedge\tau_{\varepsilon,M}^{n}}
(|v_{n}^{\varepsilon}(s)|^{2}+||v_{n}^{\varepsilon}(s)||^{6})ds\notag\\
&+2\varepsilon C\int_{0}^{t\wedge\tau_{\varepsilon,M}^{n}}|v^{\varepsilon}_{n}(s)|^{\beta+1}_{\beta+1}ds\notag\\
&+2\sqrt{\varepsilon}|\int_{0}^{t\wedge\tau_{\varepsilon,M}^{n}}(u^{\varepsilon}_{n}(s)-v_{n}^{\varepsilon}(s),
(G(\varepsilon s,u^{\varepsilon}_{n}(s))-G(\varepsilon s,v^{\varepsilon}_{n}(s)))dW(s))|)\notag\\
&\times e^{c\int_{0}^{t\wedge\tau^{n}_{\varepsilon,M}}||u^{\varepsilon}_{n}(t)||^{4}ds+\varepsilon Kt}.
\end{align}
Similarly to the proof of (5.24) and it follows from (5.41) that
\begin{align}
&(\mathbb{E}\sup_{0\leq s\leq 1\wedge\tau_{\varepsilon,M}^{n}}|u^{\varepsilon}_{n}(s)-v_{n}^{\varepsilon}(s)|^{2p})^{\frac{2}{p}}\notag\\
&\leq e^{(2cM^{2}+2\varepsilon K)} \cdot(2c^{2}\varepsilon^{2}(\mathbb{E}(\int_{0}^{t\wedge\tau^{n}_{\varepsilon,M}}| v^{\varepsilon}_{n}(s)|^{2}ds)^{p})^{\frac{2}{p}}+2c^{2}\varepsilon^{2}(\mathbb{E}(\int_{0}^{t\wedge\tau_{\varepsilon,M}^{n}}||v^{\varepsilon}_{n}(s)
||^{6}ds)^{p})^{\frac{2}{p}}\notag\\
&+8C^{2}\varepsilon^{2}(\mathbb{E}(\int_{0}^{t\wedge\tau_{\varepsilon,M}^{n}}
| v^{\varepsilon}_{n}(s)|^{\beta+1}_{\beta+1}ds)^{p})^{\frac{2}{p}}+32\varepsilon^{2}(\mathbb{E}(\int_{0}^{t\wedge\tau^{n}_{\varepsilon,M}}|| v^{\varepsilon}_{n}(s)||^{2}ds)^{p})^{\frac{2}{p}}\notag\\
&+8C\varepsilon pK^{2}\int_{0}^{t}(\mathbb{E}\sup_{0\leq r\leq s\wedge\tau_{\varepsilon,M}^{n}}|u_{n}^{\varepsilon}(r)-v_{n}^{\varepsilon}(r)|^{2p})^{\frac{2}{p}}ds)\notag\\
&\leq e^{(2cM^{2}+2\varepsilon K)}\cdot(2c^{2}M^{2}\varepsilon^{2}+2c^{2}M^{6}\varepsilon^{2}+8C^{2}\varepsilon^{2} M^{2}+32\varepsilon^{2}M^{2})\notag\\
&+e^{(2cM^{2}+2\varepsilon K)}\cdot8C\varepsilon pK^{2}\int_{0}^{t}(\mathbb{E}\sup_{0\leq r\leq s\wedge\tau_{\varepsilon,M}^{n}}|u_{n}^{\varepsilon}(r)-v_{n}^{\varepsilon}(r)|^{2p})^{\frac{2}{p}}ds.
\end{align}
Let $C_{\varepsilon,M}=e^{(2cM^{2}+2\varepsilon K)}\cdot8C\varepsilon pK^{2}$, we get
\begin{align}
&(\mathbb{E}\sup_{0\leq s\leq 1\wedge\tau_{\varepsilon,M}^{n}}|u^{\varepsilon}_{n}(s)-v_{n}^{\varepsilon}(s)|^{2p})^{\frac{2}{p}}\notag\\
&\leq e^{(2cM^{2}+2\varepsilon K)}\cdot(2c^{2}M^{2}\varepsilon^{2}+2c^{2}M^{6}\varepsilon^{2}+8C^{2}\varepsilon^{2} M^{2}+32\varepsilon^{2}M^{2})\cdot e^{C_{\varepsilon,M}}.
\end{align}
From (5.17) and Lemma 5.1 and Lemma 5.2, for any $R>0$, there exists a positive constant $M$ such that
\begin{eqnarray}
\sup_{0<\varepsilon\leq 1}\varepsilon log P((|u^{\varepsilon}_{n}|^{H}_{V}(1))^{2}> M)\leq-R,\\
\sup_{0<\varepsilon\leq 1}\varepsilon log P(\sup_{0<\varepsilon\leq 1}||u_{n}^{\varepsilon}||^{2}> M)\leq-R,\\
\sup_{0<\varepsilon\leq 1}\varepsilon log P(\sup_{0<\varepsilon\leq 1}||v_{n}^{\varepsilon}||^{2}> M)\leq-R,\\
\sup_{0<\varepsilon\leq 1}\varepsilon log P(\sup_{0<\varepsilon\leq 1}| v_{n}^{\varepsilon}|^{\beta+1}_{\beta+1}> M)\leq-R.
\end{eqnarray}
Let $p=\frac{2}{\varepsilon}$ in (5.42), we have the following inequality
\begin{align}
&\varepsilon logP(\sup_{0\leq t\leq1}|u_{n}^{\varepsilon}(t)-v_{n}^{\varepsilon}(t)|^{2}>\delta,(|u_{n}^{\varepsilon}|_{V}^{H}(1))^{2}\leq M,\sup_{0\leq t\leq1}||u_{n}^{\varepsilon}(t)||^{2}\leq M,\notag\\
&\sup_{0\leq t\leq1}||v_{n}^{\varepsilon}(t)||^{2}\leq M,\sup_{0\leq t\leq1}| v_{n}^{\varepsilon}(t)|^{\beta+1}_{\beta+1}\leq M)\notag\\
&\leq \varepsilon logP(\sup_{0\leq t\leq1\wedge\tau^{n}_{\varepsilon,M}}|u_{n}^{\varepsilon}(t)-v_{n}^{\varepsilon}(t)|^{2}>\delta)\notag\\
&\leq  log(\mathbb{E}\sup_{0\leq t\leq1\wedge\tau^{n}_{\varepsilon,M}}|u_{n}^{\varepsilon}(t)-v_{n}^{\varepsilon}(t)|^{2p})^{\frac{2}{p}}-log\delta^{2}\notag\\
&\leq 2cM^{2}+2\varepsilon K+log(2c^{2}M^{2}\varepsilon^{2}+2c^{2}M^{6}\varepsilon^{2}+8C^{2}\varepsilon^{2} M^{2}+32\varepsilon^{2}M^{2})\notag\\
&+C_{\varepsilon,M}-log\delta^{2}\notag\\
&\rightarrow -\infty,~~as~\varepsilon\rightarrow0.
\end{align}
For any $M$, (5.48) implies that there exists a positive number $\varepsilon_{0}$, such that for any $\varepsilon<\varepsilon_{0}$,
\begin{align}
&P(\sup_{0\leq t\leq1}|u_{n}^{\varepsilon}(t)-v_{n}^{\varepsilon}(t)|^{2}>\delta,(|u_{n}^{\varepsilon}|_{V}^{H}(1))^{2}\leq M,\sup_{0\leq t\leq1}||u_{n}^{\varepsilon}(t)||^{2}\leq M,\notag\\
&\sup_{0\leq t\leq1}||v_{n}^{\varepsilon}(t)||^{2}\leq M,\sup_{0\leq t\leq1}| v_{n}^{\varepsilon}(t)|^{\beta+1}_{\beta+1}\leq M)\leq e^{-\frac{R}{\varepsilon}}.
\end{align}
Combining (5.44), (5.45), (5.46), (5.47) and (5.49) together, we can get a positive number $\varepsilon_{0}$, such that for $\varepsilon\leq\varepsilon_{0}$,
\begin{eqnarray*}
P(\sup_{0\leq t\leq1}|u_{n}^{\varepsilon}(t)-v_{n}^{\varepsilon}(t)|^{2}>\delta)\leq 3e^{-\frac{R}{\varepsilon}}.
\end{eqnarray*}
Since $R$ is arbitrary, this completes the proof of Lemma 5.5.\\
{\bf Proof of (5.6)} By Lemma 5.3 and Lemma 5.4, we get for any $R>0$ that there exists an $N_{0}$ satisfying the following inequality
\begin{eqnarray}
P(\sup_{0\leq t\leq1}|u^{\varepsilon}(t)-u_{N_{0}}^{\varepsilon}(t)|^{2}>\delta)\leq e^{-\frac{R}{\varepsilon}},~~for~any~\varepsilon\in(0,1],
\end{eqnarray}
\begin{eqnarray}
P(\sup_{0\leq t\leq1}|v^{\varepsilon}(t)-v_{N_{0}}^{\varepsilon}(t)|^{2}>\delta)\leq e^{-\frac{R}{\varepsilon}},~~for~any~\varepsilon\in(0,1].
\end{eqnarray}
By virtue of Lemma 5.5, for any $N_{0}$, there exists $\varepsilon_{0}$ such that for any $\varepsilon\in(0,\varepsilon_{0}]$,
\begin{eqnarray}
P(\sup_{0\leq t\leq1}|u_{N_{0}}^{\varepsilon}(t)-v_{N_{0}}^{\varepsilon}(t)|^{2}>\delta)\leq e^{-\frac{R}{\varepsilon}}.
\end{eqnarray}
Hence, for any $\varepsilon\in(0,\varepsilon_{0}]$,
\begin{eqnarray}
P(\sup_{0\leq t\leq1}|u^{\varepsilon}(t)-v^{\varepsilon}(t)|^{2}>\delta)\leq 3e^{-\frac{R}{\varepsilon}}.
\end{eqnarray}
Since $R$ is arbitrary, we deduce
\begin{eqnarray*}
\lim_{\varepsilon\rightarrow0}\varepsilon logP(\sup_{0\leq t\leq1}|u^{\varepsilon}(t)-v^{\varepsilon}(t)|^{2}>\delta)=-\infty.
\end{eqnarray*}

\end{document}